\input  amstex
\input amsppt.sty
\magnification1200
\vsize=23.5truecm
\hsize=16.5truecm

\NoBlackBoxes

\def\supp{\operatorname{supp}}

\def\crm{\overline{\Bbb R}_-}

\def\rnp{{\Bbb R}^n_+}

\def\rnpm{\Bbb R^n_\pm}
\def\crnp{\overline{\Bbb R}^n_+}

\def\crnpm{\overline{\Bbb R}^n_\pm}
\def\comega{\overline\Omega }
\def\ang#1{\langle {#1} \rangle}
\def\simto{\overset\sim\to\rightarrow}

\def\rp{ \Bbb R_+}

\def\R{\Bbb R}

\def\ol{\overline}
\def\SD{\Cal S}
\def\F{\Cal F}
\def\E{\Cal E}
\def\D{\Cal D}

\document

\medskip

\topmatter
\title
Local and nonlocal boundary conditions for $\mu $-transmission and
fractional  elliptic pseudodifferential operators
\endtitle
\author Gerd Grubb \endauthor 
\affil
{Department of Mathematical Sciences, Copenhagen University,
Universitetsparken 5, DK-2100 Copenhagen, Denmark.}
E-mail {\tt grubb\@math.ku.dk}
\endaffil
\rightheadtext{Boundary problems}
\abstract
A classical pseudodifferential operator $P$ on ${\Bbb R}^n$
satisfies the $\mu
$-transmission condition relative to a smooth open subset $\Omega $,
 when the symbol terms have a certain twisted
 parity on the normal to $\partial\Omega
$. As shown recently by the author, the condition assures 
solvability of Dirichlet-type
boundary problems for $P$ in full scales of Sobolev spaces with a singularity
$d^{\mu -k}$, $d(x)=\operatorname{dist}(x,\partial\Omega )$. Examples include
fractional Laplacians $(-\Delta )^a$ and
complex powers of strongly elliptic PDE. 

We now introduce new boundary conditions, of Neumann type or
more general nonlocal. It is also shown how problems with data on ${\Bbb
R}^n\setminus \Omega $ reduce to problems supported on
$\comega$, and how the so-called ``large'' solutions arise. Moreover,
the results are extended to general function
spaces $F^s_{p,q}$ and $B^s_{p,q}$, including  H\"older-Zygmund spaces
$B^s_{\infty ,\infty }$. This leads to optimal
H\"older estimates, e.g.\ for Dirichlet solutions of $(-\Delta
)^au=f\in L_\infty (\Omega )$,
 $u\in d^aC^a(\comega)$ when $0<a<1$, $a\ne\frac12$.

\endabstract

\endtopmatter

Boundary value problems for elliptic pseudodifferential operators
($\psi $do's) $P$, on a smooth subset $\Omega $ of a Riemanninan
manifold $\Omega _1$, have been studied under various hypotheses through the
years. There is a well-known calculus initiated by Boutet de Monvel \cite{B71, RS82,
G84, G90,  G96, S01, G09} for integer-order $\psi $do's with the 0-transmission
property (preserving $C^\infty $ up to the boundary), including
boundary value problems for
elliptic differential operators and their inverses. There are theories
treating more general
operators with suitable factorizations of the
principal symbol, initiated by Vishik and Eskin, see e.g.\ \cite{E81, S94, CD01}. 
Theories for operators without the transmission property have been developed by Schulze and coauthors, see e.g.\ \cite{RS84, HS08}, and
theories where the boundary is considered as a singularity of the
manifold have been developed in works of Melrose and coauthors, see e.g. \cite{M93,
AM09}. 

A category of $\psi $do's lying between the operators handled by the
Boutet de Monvel calculus and the very general categories mentioned
above, consists of the $\psi $do's with a $\mu $-transmission
property, $\mu \in{\Bbb C}$, with respect to $\partial\Omega $. Only
recently, a systematic study in $H^s_p$ Sobolev spaces was given in Grubb \cite{G15a},
departing from a result on such operators in $C^\infty $-spaces by
H\"ormander \cite{H85} Th.\ 18.2.18 (in fact developed from a lecture
note of H\"ormander \cite{H65}). This category includes fractional
Laplacians $(-\Delta )^a$ and complex powers of strongly elliptic
differential operators, and also more generally polyhomogeneous $\psi
$do's with symbol $p\sim\sum_{j\in{\Bbb N}_0}p_j$ having {\it even parity}
($p_j(x,-\xi )=(-1)^jp_j(x,\xi )$ for $j\ge 0$) or a twisted parity involving a factor
$e^{i\pi \varrho }$.  The general $\mu $-transmission operators have such a
reflection property of the symbol at $\partial\Omega $ just in the
normal direction, see (1.5) below. It allows regularity and
solvability results not only for $s$ in an interval, but for all $s\to\infty $.

The fractional Laplacian and its generalizations, often formulated as
singular integral operators, 
are currently of interest 
both in probability theory
and finance, in mathematical physics and in geometry.

\medskip
The work \cite{G15a} showed the Fredholm solvability of homogeneous or
nonhomogeneous Dirichlet-type
problems in large scales of Sobolev spaces, for $\mu $-transmission $\psi
$do's. In the present paper we introduce more general boundary
conditions and find criteria for their solvability. There are the general
nonlocal conditions $\gamma _0Bu=\psi $, where $B$ is a $\mu
$-transmission $\psi $do; in addition to this, local higher-order
conditions such as a Neumann-type condition involving the normal
derivative at $\partial\Omega $ are treated. The case of
matrix-formed $P$ is briefly considered.

Moreover, we show by use of Johnsen \cite{J96} that the theory also
works in Besov-Triebel-Lizorkin spaces $B^s_{p,q}$ and $F^s_{p,q}$,
with special attention to the spaces $B^s_{\infty ,\infty }$ that
coincide with 
H\"older spaces $C^s$ for $s\in \rp\setminus {\Bbb N}$. In comparison
with \cite{G15a}, this
allows for a sharpening of H\"older results for $(-\Delta )^a$ (and other
$a$-transmission operators), as follows:
Let $\comega$ be compact $\subset{\Bbb R}^n$. For
solutions $u\in e^+L_\infty (\Omega )$ of $r^+(-\Delta )^au=f$,
$$
f\in L_\infty (\Omega )\implies u\in e^+d(x)^aC^a(\comega), \text{
when }a\in \,]0,1[\, , \; a\ne \tfrac12, \tag0.1
$$
which is optimal in the H\"older exponent. (For $a=\frac12$, it holds
with $C^a$ replaced by $C^{a-\varepsilon }$. Also higher regularities
are treated, and optimal H\"older estimates for nonhomogeneous
Dirichlet and Neumann problems are likewise shown.) In a new work \cite{RS14},
Ros-Oton and Serra have studied integral operators with homogeneous, positive,
even kernel and obtained (0.1) with $C^a$ replaced by
$C^{a-\varepsilon }$; in the smooth case this is covered by the
present theory. (We are concerned with linear operators; the nonlinear implications in
\cite{RS14}  are not touched here.) Such operators 
were treated in cases without boundary by Caffarelli and Silvestre,
see e.g.\ \cite{CS09}.

Furthermore, we show the equivalence of Dirichlet problems for $u$
supported in $\comega $ with problems
prescribing a value of $u$ on the exterior ${\Bbb R}^n\setminus\Omega
$, obtaining new results for the latter, that were
treated recently by e.g.\  Felsinger, Kassman and Voigt in \cite{FKV14} and
Abatangelo in \cite{A14}. 

For nonhomogeneous problems the solutions can be ``large'' at the
boundary, cf.\ \cite{A14} and its references. We show how the
solutions have a specific power singularity when the
boundary data are nontrivial.

The case $a=\frac12$ enters as a boundary integral operator in treatments of mixed boundary value problems
for elliptic differential operators. 
The present results are applied to mixed problems in \cite{G15b}.
\medskip

{\it Outline.} In Section 1, we recall briefly the relevant
definitions of operators and spaces. Section 2 presents the basic results
on Dirichlet and Neumann problems for $(-\Delta )^a$, including
situations with given exterior data, and deriving conclusions in
H\"older spaces. Section 3 explains
the extension of the general results to Besov-Triebel-Lizorkin spaces,
including $B^s_{\infty ,\infty }$. Section 4 introduces new nonlocal
boundary conditions $\gamma _0Bu=\psi $, as well as local Neumann-type
conditions; also matrix-formed $P$ are discussed. The Appendix illustrates the theory by treating a particular constant-coefficient
case, showing how the problems for $(1-\Delta
)^a$ on $\rnp$ can be solved in full detail by explicit calculations.

\head 1. Preliminaries
\endhead

The notations of \cite{G15a} will be used. We shall give a brief
account here, and refer there for further details.

Consider a Riemannian $n$-dimensional $C^\infty $ manifold $\Omega
_1$ (it can be ${\Bbb R}^n$) and an embedded smooth $n$-dimensional
manifold $\comega$ with boundary $\partial\Omega $ and interior
$\Omega $. For $\Omega _1={\Bbb R}^n$,  $\Omega $ can be $\rnpm=\{x\in
{\Bbb R}^n\mid x_n\gtrless 0\}$; here $(x_1,\dots, x_{n-1})=x'$. In
the general manifold case, $\comega$ is taken compact. For
$\xi \in{\Bbb R}^n$, we
denote $(1+|\xi |^2)^\frac12 =\ang \xi $, and denote by $[\xi ]$ a positive
$C^\infty $-function equal to $|\xi |$ for $|\xi |\ge 1$ and $\ge
\tfrac12$ for all $\xi $. Restriction from $\R^n$ to $\rnpm$ (or from
$\Omega _1$ to $\Omega $ resp.\ $\complement\comega$) is denoted $r^\pm$,
 extension by zero from $\rnpm$ to $\R^n$ (or from $\Omega $ resp.\
 $\complement\comega$ to $\Omega _1$) is denoted $e^\pm$.

A pseudodifferential operator ($\psi $do) $P$ on ${\Bbb R}^n$ is
defined from a symbol $p(x,\xi )$ on ${\Bbb R}^n\times{\Bbb R}^n$ by 
$$
Pu=p(x,D)u=\operatorname{OP}(p(x,\xi ))u 
=(2\pi )^{-n}\int e^{ix\cdot\xi
}p(x,\xi )\hat u\, d\xi =\Cal F^{-1}_{\xi \to x}(p(x,\xi )\hat u(\xi ));\tag1.1
$$  
here $\Cal F$ is the Fourier transform $(\F u)(\xi )=\hat u(\xi
)=\int_{{\Bbb R}^n}e^{-ix\cdot \xi }u(x)\, dx$. The symbol $p$ is
assumed to be such that $\partial_x^\beta \partial_\xi ^\alpha p(x,\xi
)$ is $O(\ang\xi ^{r-|\alpha |})$ for all $\alpha ,\beta $, for some
$r\in{\Bbb R}$ (defining
the symbol class $S^r_{1,0}({\Bbb R}^n\times{\Bbb R}^n)$); then it has
order $r$. The definition of $P$
is carried over to manifolds by use of local coordinates. We refer to
textbooks such as \cite{H85}, Taylor \cite{T91}, \cite{G09} for the rules of
calculus; \cite{G09} moreover gives an account of the Boutet de Monvel
calculus of pseudodifferential boundary problems, cf.\ also e.g.\
\cite{G96}, Schrohe \cite{S01}.  When $P$ is a $\psi $do on ${\Bbb R}^n$ or $\Omega _1$, $P_+=r^+Pe^+$
denotes its truncation to $\rnp$ resp.\ $\Omega $.
\medskip

Let $1<p<\infty $  (with $1/p'=1-1/p$), then we define for $s\in{\Bbb R}$ the spaces
$$
\aligned
H^s_p(\R^n)&=\{u\in \SD'({\Bbb R}^n)\mid \F^{-1}(\ang{\xi }^s\hat u)\in
L_p(\R^n)\},\\
\dot H^{s}_p(\crnp)&=\{u\in H^{s}_p({\Bbb R}^n)\mid \supp u\subset
\crnp \},\\
\ol H^{s}_p(\rnp)&=\{u\in \D'(\rnp)\mid u=r^+U \text{ for some }U\in
H^{s}_p(\R^n)\};
\endaligned \tag1.2
$$
here $\operatorname{supp}u$ denotes the support of $u$. For $\comega$
compact $\subset\Omega _1$, the definition extends to define $\dot
H^s_p(\comega)$ and $\ol H^s_p(\Omega )$ by use of a finite system of
local coordinates. We shall in the present paper moreover work in the
Triebel-Lizorkin and Besov  spaces $F^s_{p,q}$ and $B^s_{p,q}$, defined for
$s\in{\Bbb R}$, $0<p,q\le \infty $ (we take $p<\infty $ in the $F$-case), and
the derived spaces $\dot F^s_{p,q}$ and $\ol F^s_{p,q}$, etc. Here we
refer to Triebel \cite{T95} and Johnsen \cite{J96} for the basic
explanations. (\cite{T95} writes $\widetilde F$ instead of $\dot
F$, etc., the present notation stems from H\"ormander's works.)
For
H\"older spaces $C^t$, $\dot C^t(\comega)$ denotes the H\"older
function on $\Omega _1$ supported in $\comega$. 
$B^s_{p,p}$ is also denoted $B^s_p$ when $p<\infty $, and $F^s_{p,p}=B^s_{p,p}$,
$F^s_{p,2}=H^s_p$, $H^s_2=B^s_2$.

We shall use the conventions  $\bigcup_{\varepsilon >0}
H_p^{s+\varepsilon } =
H_p^{s+ 0}$, $\bigcap_{\varepsilon >0} H_p^{s-\varepsilon }=
H_p^{s- 0}$, applied in a similar way for the other scales of
spaces.

The results hold in particular for 
$B^s_{\infty ,\infty }$-spaces. These are interesting because $B^s_{\infty
,\infty }({\Bbb R}^n)$ equals
the H\"older space $C^s({\Bbb R}^n)$ when $s\in\rp\setminus {\Bbb
N}$. (There are similar
statements for derived spaces over $\rnp$ and $\Omega $.) 
The spaces $B^s_{\infty ,\infty }({\Bbb R}^n)$ identify with the
H\"older-Zygmund spaces, often denoted $\Cal C^s({\Bbb R}^n)$ when
$s>0$.
There is a nice account of these spaces in Section 8.6 of \cite{H97},
where they are denoted $C^s_*({\Bbb R}^n)$ for all $s\in{\Bbb R}$; we
shall use this label below, for simplicity of notation:
$$B^s_{\infty ,\infty }=C^s_*,\quad \text{ for
all }s\in{\Bbb R}.\tag 1.3$$ 
For integer values one has, with $C^k_b({\Bbb R}^n)$ denoting the space of functions with bounded
continuous derivatives up to order $k$,
$$
\aligned
C^k_b({\Bbb R}^n)&\subset C^{k-1,1}({\Bbb R}^n)\subset C^k_*({\Bbb R}^n)\subset C^{k-0}({\Bbb R}^n)\text{ when }k\in
{\Bbb N},\\
 \quad C^0_b({\Bbb R}^n)&\subset L_\infty ({\Bbb R}^n)\subset
 C^0_*({\Bbb R}^n),
\endaligned\tag1.4
$$
and similar statements for derived spaces.
\medskip

A $\psi $do $P$ is called classical (or polyhomogeneous) when the
symbol $p$ has an asymptotic expansion $p(x,\xi )\sim \sum_{j\in{\Bbb
N}_0}p_j(x,\xi )$ with $p_j$ homogeneous in $\xi $ of degree $m-j$ for
all $j$. Then $P$ has order $m$. One can even allow $m$ to be complex;
then $p\in S^{\operatorname{Re} m}_{1,0}({\Bbb R}^n\times{\Bbb R}^n)$;
the operator and symbol are still said to be of order $m$. 

Here there is an additional definition: $P$ satisfies {\it the $\mu
$-transmission condition} (in short: is of type $\mu $) for some $\mu
\in{\Bbb C}$  when, in local coordinates,
$$
\partial_x^\beta \partial_\xi ^\alpha p_j(x,-N)=e^{\pi i(m-2\mu -j-|\alpha | )
}\partial_x^\beta \partial_\xi ^\alpha p_j(x,N),\tag1.5
$$
for all $x\in\partial\Omega $, all $j,\alpha ,\beta $, where 
$N$ denotes the interior normal to $\partial\Omega $ at $x$.
The implications of the $\mu $-transmission property 
were a main subject of \cite{G15a}.
  
A special role in the theory is played by the {\it order-reducing
operators}. There is a simple definition of operators $\Xi _\pm^\mu $ on
${\Bbb R}^n$
$$ 
\Xi _\pm^\mu =\operatorname{OP}(([\xi ']\pm i\xi _n)^\mu ) 
$$
(or with $[\xi ']$ replaced by $\ang{\xi '}$); they preserve support
in $\crnpm$, respectively. Here the function
$([\xi ']\pm i\xi _n)^\mu $ does not satisfy all the estimates
required for the class $S^{\operatorname{Re}\mu }({\Bbb
R}^n\times{\Bbb R}^n)$, but the operators are useful for some
purposes. There is a more refined choice $\Lambda _\pm^\mu $ (with
symbol $\lambda _\pm^\mu (\xi )$) that does
satisfy all the estimates, and there is a definition $\Lambda
_\pm^{(\mu )}$ in the manifold situation. These operators 
 define
homeomorphisms for all $s\in{\Bbb R}$ such as 
$$
\aligned
\Lambda^{(\mu )}_+&\colon \dot H^s_p(\comega )\simto
\dot H^{s-\operatorname{Re}\mu }_p(\comega ),\\
\Lambda ^{(\mu )}_{-,+}&\colon \ol H^s_p(\Omega )\simto
\ol H^{s-\operatorname{Re}\mu }_p(\Omega );
\endaligned \tag1.6
$$
 here $\Lambda ^{(\mu )}_{-,+}$ is short for $r^+\Lambda ^{(\mu
 )}_{-}e^+$, suitably extended to large negative $s$ (cf.\ Rem.\ 1.1
 and Th.\ 1.3 in \cite{G15a}).

The following special spaces introduced by H\"ormander are particularly adapted to $\mu
$-transmission operators $P$: 
$$
\aligned
H^{\mu (s)}_p(\crnp)&=\Xi _+^{-\mu }e^+\ol H_p^{s-\operatorname{Re}\mu
}(\rnp),\quad  s>\operatorname{Re}\mu -1/p',\\
H^{\mu (s)}_p(\comega)&=\Lambda  _+^{(-\mu )}e^+\ol H_p^{s-\operatorname{Re}\mu
}(\Omega ),\quad  s>\operatorname{Re}\mu -1/p',\\
\Cal E_\mu (\comega)&=e^+\{u(x)=d(x)^\mu v(x)\mid v\in C^\infty
(\comega)\};
\endaligned\tag1.7
$$
namely, $r^+P$ (of order $m$) maps them into $\ol H_p^{s-\operatorname{Re}m}(\rnp)$, $\ol
 H_p^{s-\operatorname{Re}m}(\Omega)$ resp.\ $C^\infty (\comega)$ (cf.\
 \cite{G15a} Sections 1.3, 2, 4), and they appear as domains of
 realizations of $P$ in the elliptic case. In the third line, $\operatorname{Re}\mu
 >-1$ (for other $\mu $, cf.\ \cite{G15a}) and $d(x)$ is a $C^\infty $-function vanishing to
 order 1 at $\partial\Omega $ and positive on $\Omega $, e.g.\
 $d(x)=\operatorname{dist}(x,\partial\Omega )$ near $\partial\Omega $.
 One has that $H^{\mu
 (s)}_p(\comega)\supset \dot H_p^s(\comega)$, and the distributions are
 locally in $H^s_p$   on $\Omega $, but at the boundary they in general have a
 singular behavior. More about that in the text below.

The order-reducing operators also operate in the Besov-Triebel-Lizorkin
scales of spaces, satisfying the relevant versions of (1.6), and the definitions
in (1.7) extend.

\head 2. Three basic problems for the fractional Laplacian
\endhead

As a useful introduction, we start out by giving a detailed presentation of
boundary problems for the basic example of
the fractional Laplacian. 

Let $P_a=(-\Delta )^a$, $a>0$, and let $\Omega $ be a bounded open
subset of ${\Bbb R}^n$ with a $C^\infty $-boundary $\partial\Omega =\Sigma $. $P_a$, acting as
$u\mapsto \F^{-1}(|\xi |^{2a}\hat u)$, is a pseudodifferential
operator on ${\Bbb R}^n$ of order $2a$, and it is of type $a$ and has
factorization index $a$ relative to $\Omega $, as defined in
\cite{G15a}. With a terminology
introduced by H\"ormander in the notes \cite{H65} and now exposed in
\cite{G15a}, we consider the following problems for $P_a$:
\roster
\item {\it The homogeneous Dirichlet problem}
$$
\cases r^+P_au& =f\text{ on }\Omega ,\\
\operatorname{supp}u&\subset \comega. \endcases \tag2.1
$$
\item {\it A nonhomogeneous Dirichlet problem} (with $u$ less regular
than in (2.1))
$$
\cases r^+P_au& =f\text{ on }\Omega ,\\
\operatorname{supp}u&\subset \comega ,\\
d(x)^{1-a}u &=\varphi \text{ on }\Sigma .
\endcases \tag2.2
$$
\item {\it A nonhomogeneous Neumann problem}
$$
\cases r^+P_au& =f\text{ on }\Omega ,\\
\operatorname{supp}u&\subset \comega ,\\
\partial_n(d(x)^{1-a}u) &=\psi  \text{ on }\Sigma  .
\endcases \tag2.3
$$
\endroster
It is shown in \cite{G15a} that (2.1) and (2.2) have good 
solvability properties in suitable Sobolev spaces and H\"older spaces,
and we shall include (2.3) in the study below. In the following, we derive
further properties of each of the three problems. 

\example{Remark 2.1} The theorems in Sections 2.1 and 2.2 below are also valid when $(-\Delta
)^a$ is replaced by a general
$a$-transmission $\psi $do $P$ of order $2a$ and with factorization
index $a$, except that the bijectiveness is replaced by the Fredholm
property. They also hold when $\comega$ is a compact subset of a manifold $\Omega _1$.  The results in Section 2.3 extend to such operators when
they are principally like $(-\Delta )^a$. 
\endexample

In the appendix of this paper we have included a treatment of
$(1-\Delta )^a$ on a half-space; it is a model case where one can
obtain the solvability results directly by Fourier transformation.

\subhead 2.1 The homogeneous Dirichlet problem \endsubhead

From a point of view of functional analysis (as used e.g.\ in Frank
and Geisinger \cite{FG14}), it is natural to define
the Dirichlet realization $P_{a,D}$ as the Friedrichs extension of the
symmetric operator $P_{a,0}$ in $L_2(\Omega )$ acting like $r^+P_a$
with domain $C_0^\infty (\Omega )$. There is an associated
sesquilinear form 
$$
p_{a,0}(u,v)=(2\pi )^{-n}\int_{{\Bbb R}^n}|\xi |^{2a}\hat u(\xi
)\bar{\hat v(\xi )}\,d\xi ,\quad u,v\in C_0^\infty (\Omega ).\tag2.4
$$
Since $(\|u\|^2_{L_2}+\int |\xi |^{2a}|\hat u|^2\,d\xi )^\frac12$ is a
norm equivalent with $\|u\|_{H_2^{a}}$, the completion of $C_0^\infty
(\Omega )$ in this norm is $V=\dot H_2^a(\comega)$, and $p_{a,0}$ extends to
a continuous nonnegative symmetric sesquilinear form on $V$. A standard
application of the Lax-Milgram lemma (e.g.\ as in \cite{G09}, Ch.\ 12)
gives the operator $P_{a,D}$ that
is selfadjoint nonnegative in $L_2(\Omega )$ and acts like
$r^+P_a:\dot H_2^a(\comega)\to \ol H_2^{\,-a}(\Omega )$, with domain 
$$
D(P_{a,D})=\{u\in \dot H_2^a(\comega)\mid r^+P_au\in L_2(\Omega )\}.\tag2.5
$$
The operator has compact resolvent, and the spectrum is a
nondecreasing sequence of nonnegative eigenvalues going to
infinity. As we shall document below, 0 is not an eigenvalue, so
$P_{a,D}$ in fact has a positive lower bound and is invertible. 

The results of \cite{G15a} (Sections 4, 7) 
clarify the mapping
properties and solvability properties further: For $1<p<\infty $,  $r^+P_a$ maps continuously
$$
r^+P_a\colon H_p^{a(s)}(\comega)\to \ol H_p^{s-2a}(\Omega ), \text{
when }s>a-1/p';\tag2.6
$$
there is the regularity result 
$$
u\in \dot H_p^{a-1/p'+0} (\comega),\; r^+ P_au\in \ol
H_p^{s-2a}(\Omega )\implies u\in H_p^{a(s)}(\comega),\text{ when
}s>a-1/p' ,
\tag2.7 
$$
and the mapping (2.6)  is Fredholm. (It is even bijective, as seen
below.) As an application of the results for $s=2a$,
$p=2$, we have in particular that
$$
D(P_{a,D})=H_2^{a(2a)}(\comega)=\Lambda _+^{(-a)}e^+\ol H_2^{a}(\Omega ),\tag2.8
$$
see also Example 7.2 in \cite{G15a}.
We recall from \cite{G15a} Th.\ 5.4 that
$$
H_p^{a(s)}(\comega)\cases =\dot H_p^s(\comega),\text{ when
}a-1/p'<s<a+1/p,\\
\subset\dot H_p^{s-0}(\comega),\text{ when }s=a+1/p,\\
\subset e^+d^a\ol H_p^{s-a}(\Omega)+ \dot H_p^{s}(\comega),\text{ when
}s>a+1/p, s-a-1/p\notin {\Bbb N},\\
\subset e^+d^a\ol H_p^{s-a}(\Omega)+ \dot H_p^{s-0}(\comega),\text{ when
}s-a-1/p\in{\Bbb N}.
\endcases \tag2.9
$$

In \cite{G15a} we used Sobolev embedding theorems to draw conclusions for
H\"older spaces, cf.\ Section 7 there. Slightly sharper (often
optimal) results can be
obtained if we go via an extension of the results of \cite{G15a} to the
general scales of Triebel-Lizorkin and Besov spaces
$F^s_{p,q}$ and $B^s_{p,q}$.  The
extended theory will be presented in detail 
below in Sections 3--4; for the moment
we shall borrow some results to give powerful statements for $(-\Delta )^a$,
$0<a<1$. We recall that the notation $B^s_{\infty ,\infty }$ is
simplified to $C^s_*$, and that
$C^s_*$ equals $C^s$ (the ordinary H\"older space) for  $s\in\rp\setminus {\Bbb N}$, cf.\
also (1.4). Moreover, as special cases of Definition 3.1 and Theorem
3.4 below for $p=q=\infty $,
$$
\aligned
C^{\mu  (s)}_*(\comega)&=\Lambda _+^{(-\mu 
)}e^+\ol C^{s-\operatorname{Re}\mu  }_*(\Omega )\text{ for
}s>\operatorname{Re}\mu  -1,\text{ and}\\
C^{\mu  (s)}_*(\comega)&\subset \cases 
d(x)^\mu 
e^+\ol C^{s-\operatorname{Re}\mu  }_*(\Omega )+ \dot
C^s_*(\comega)\text{ when }s>\operatorname{Re}\mu  ,\, s-\operatorname{Re}\mu \notin{\Bbb
N},\\
d(x)^\mu 
e^+\ol C^{s-\operatorname{Re}\mu  }_*(\Omega )+ \dot
C^{s-0}_*(\comega)\text{ when }s>\operatorname{Re}\mu  ,\, s-\operatorname{Re}\mu \in{\Bbb N}.\endcases
\endaligned\tag2.10
$$
Note also that the distributions in  $C^{\mu  (s)}_*(\comega)$ are locally in $C^s_*
$ on $\Omega $, by the ellipticity of $\Lambda _+^{(-\mu  )}$.

{\it We focus in the following on the case $0<a<1$, assumed in the
rest of this chapter.} Here we find the following results, with conclusions formulated in ordinary H\"older spaces:

\proclaim{Theorem 2.2} Let $s>a-1$.  If $u\in \dot C^{a-1+\varepsilon
}_*(\comega)$ for some $\varepsilon >0$ (e.g.\ if $u\in e^+L_\infty (\Omega )$), and
$r^+Pu\in \ol
C_*^{s-2a}(\Omega)$, then $u\in C^{a(s)}_*(\comega)$. The mapping $r^+P_a$ defines a bijection
$$
r^+P_a\colon C^{a(s)}_*(\comega)\to \ol
C_*^{s-2a}(\Omega).\tag2.11
$$
In particular, for any $f\in L_\infty (\Omega )$, there exists a
unique solution $u$ of {\rm (2.1)} in $C^{a(2a)}_*$; it satisfies
$$
\aligned
u&\in  e^+d(x)^aC^a(\comega)\cap C^{2a}(\Omega ),\text{ when }a\ne
\tfrac12,\\
u&\in \bigl(e^+d(x)^{\frac12}C^{\frac12}(\comega)+\dot
C^{1-0}(\comega)\bigr)\cap C^{1-0}(\Omega )\\
&\subset e^+d(x)^{\frac12}C^{\frac12-0}(\comega)\cap C^{1-0}(\Omega ),
\text{ when }a=
\tfrac12.
\endaligned
 \tag2.12
$$
 For $f\in C^t(\comega)$, $t>0$, the solution satisfies
$$
u\in \cases e^+d(x)^aC^{a+t}(\comega)\cap C^{2a+t}(\Omega ),\text{ when
}a+t\text{ and }2a+t\notin  {\Bbb N},\\
\bigl(e^+d(x)^aC^{a+t-0}(\comega)+\dot C^{2a+t-0}(\comega)\bigr)\cap C^{2a+t}(\Omega ),\text{ when
}a+t\in  {\Bbb N}\\
\bigl(e^+d(x)^aC^{a+t}(\comega)+\dot C^{2a+t-0}(\comega)\bigr)\cap C^{2a+t-0}(\Omega ),\text{ when
}2a+t\in  {\Bbb N}.
\endcases \tag2.13
$$

Also the mappings {\rm (2.6)} are bijections, for $s>a-1/p'$.
\endproclaim

\demo{Proof} The first two statements are a special case of Theorem
3.2 below (see Example 3.3), except that we have replaced the Fredholm property with
bijectiveness.
According to Ros-Oton and Serra \cite{RS12}, Prop.\ 1.1, a weak
solution (a solution in $\dot H^a_2(\comega )$) of the problem (2.1) with $f\in
L_\infty (\Omega )$ satisfies $\|u\|_{C^a}\le
C\|f\|_{L_\infty }$; in particular it is unique. For $f\in \ol H_2^{\,-a}(\Omega )$, the Fredholm
property of $r^+P_a$ from $H_2^{a(a)}(\comega)=\dot H_2^{a}(\comega )$ to $\ol H_2^{\,-a}(\Omega )$
is covered by \cite{G15a} Th.\ 7.1 with $s=a$, $p=2$. Moreover, the kernel
$\Cal N$ is
in $\E_a(\comega)$ by Theorem 3.5 below. If the kernel were nonzero,
there would exist nontrivial null-solutions $u\in \E_a(\comega)$,
contradicting the uniqueness for $f\in L_\infty (\Omega )$ mentioned
above. Thus $\Cal N=0$. Then the kernel of the Dirichlet realization $P_{a,D}$
in $L_2(\Omega )$ recalled above is likewise 0, and since it is a
selfadjoint operator with compact resolvent, it must be bijective. So
the cokernel in $L_2(\Omega )$ is
likewise 0. This shows the bijectivity of (2.6) in the case $s=2a$,
$p=2$. In view of Theorem 3.5 below, this bijectivity carries over to all the other versions,
including (2.6) for general $s>a-1/p'$, and the mapping (2.11) in
$C^s_*$-spaces for $s>a-1$.

For (2.12) we use Theorem 3.4 (as recalled in (2.10)), noting that $\ol C^a_*(\Omega
)=C^a(\comega)$, that  $\dot C^{2a}_*(\comega)=\dot C^{2a}(\comega)\subset 
d(x)^aC^a(\comega)$ when $a\ne \frac12$, and that $u\in
C^{2a}(\Omega )$ by
interior regularity when $a\ne\frac12$, with slightly weaker
statements when $a=\frac12$. 
The rest of the statements follow similarly by use of (2.10) with $\mu =a$ and the various
informations on the relation betweeen the $C^s_*$-spaces and standard H\"older spaces.
\qed\enddemo

Ros-Oton and Serra showed in \cite{RS12}, under  weaker smoothness
hypotheses on $\Omega $, the inclusion $u\in
d^aC^{\alpha }(\comega)$ for an $\alpha $ with $0<\alpha <\operatorname{min}\{a,1-a\}$, and improve it in a new work
\cite{RS14} to $\alpha =a-\varepsilon $. They observe that $\alpha >a$ cannot be
obtained, so $\alpha =a$ that we obtain in (2.12) is optimal.

We also have as shown in \cite{G15a} Th.\ 4.4, for functions $u$
supported in $\comega$ (cf.\ the first inclusion in (2.7)), 
$$
r^+P_au\in C^\infty (\comega)\iff u\in \E_a(\comega)\equiv
\{u=e^+d(x)^av(x)\mid v\in C^\infty (\comega)\}.\tag 2.14
$$
It is worth emphasizing that the functions in $\E_a$ have a nontrivially
singular behavior at $\Sigma $ when $a\notin{\Bbb N}_0 $; $e^+C^\infty
(\comega)$ and  $\E_a(\comega)$ are very different spaces. The
appearance of a factor $d^{\mu _0}$ 
when the
factorization index is $\mu _0$, was observed in $C^\infty$-situations also in \cite{E81} p.\
311 and in \cite{CD01} Th.\ 2.1.

The solution operator is
denoted $R$; its form as a composition of pseudodifferential factors
was given in \cite{G15a}.

\medskip

There is another point of view on the Dirichlet problem for $P_a$ that
we shall also discuss. In a number of papers, see e.g.\ Hoh and Jacob
\cite{HJ96}, Felsinger,
Kassman and Vogt \cite{FKV14} and their references, the Dirichlet
problem for $P_a$ (and other related operators) is formulated as
$$
\cases P_aU&=f\text{ in }\Omega ,\\
U&=g\text{ on }\complement\Omega .\endcases \tag 2.15
$$

Although the main aim is to determine $U$ on $\Omega $, the
prescription of the values of $U$ on $\complement \Omega $ is
explained as
 necessitated by the nonlocalness of $P_a$.
As observed explicitly in \cite{HJ96}, the transmission property of
Boutet de Monvel \cite{B71} is not satisfied; hence that theory of
boundary problems for pseudodifferential operators is of no help. 
But now that we have the $\mu
$-transmission calculus, it is worth investigating what the methods can give.

The case $g=0$ corresponds to the formulation (2.1). But also in general, (2.15)
can be reduced to (2.1) when the spaces are suitably chosen. For
(2.15), let  $f$ be given in $\ol H_p^{s-2a}(\Omega )$ (with $s>a-1/p'$), and let $g$ be
given in $\ol H^s_p(\complement\comega)$; then we search for $U$ in a
Sobolev space over ${\Bbb R}^n$. 

Let
$G=\ell g$ be an extension of $g$ to $H_p^s({\Bbb R}^n)$. Then $u=U-G$ must satisfy  
$$
\cases r^+P_au&=f-r^+P_aG\text{ in }\Omega ,\\
\supp u&\subset\comega.\endcases \tag 2.16
$$
Here $P_aG\in H_{p,\operatorname{loc}}^{s-2a}({\Bbb R}^n)$, so $f-r^+P_aG\in \ol
H_p^{s-2a}(\Omega )$.

According to our analysis of (2.1), there is 
a unique
solution $u=R(f-r^+P_aG)\in H_p^{a(s)}(\comega)$ of (2.16). Then
(2.15) has the solution
$U=u+G\in H_p^{a(s)}(\comega)+H^s_p({\Bbb R}^n)$. Moreover, there is
at most one solution to (2.15) in this space, for if
$U_1=u_1+G_1$ and $U_2=u_2+G_2$ are two solutions, then
$v=u_1-u_2+G_1-G_2$ is supported in $\comega$, hence lies in  $
H_p^{a(s)}(\comega)+\dot H_p^s(\comega)=H_p^{a(s)}(\comega)$ and satisfies (2.1) with $f=0$, hence it must be 0.

This reduction allows a study of higher regularity of the
solutions. The treatment in \cite{FKV14} seems primarily directed
towards the regularity involved in variational formulations ($p=2$,
$s=a$) where
Vishik and Eskin's results would be applicable; moreover, \cite{FKV14} allows a
less smooth boundary.

We have shown:

\proclaim{Theorem 2.3} Let $s>a-1/p'$, and let $f\in\ol
H_p^{s-2a}(\Omega )$ and  $g\in \ol H^s_p(\complement\comega)$ be
given. Then the problem {\rm (2.15)} has the unique solution $U=u+G\in
H_p^{a(s)}(\comega)+H^s_p({\Bbb R}^n)$, where
$G\in H^s_p({\Bbb R}^n)$ is an extension of $g$ and 
$$
u=R(f-r^+P_aG)\in H_p^{a(s)}(\comega);\tag2.17
$$
here $R$ is the solution operator for {\rm (2.1)}.
\endproclaim

Observe in particular that 
the solution is independent of the choice of an extension
operator $\ell \colon g\mapsto G$.

There is an immediate corollary for solutions in H\"older spaces (as
in \cite{G15a} Sect.\ 7):

\proclaim{Corollary 2.4} 
Let $p>n/a$. For $f\in L_p(\Omega )$, $g\in
C^{2a+0}(\complement\Omega )\cap \ol H_p^{2a}(\complement\comega )$, the
solution of {\rm (2.15)} according to Theorem {\rm 2.3} satisfies 
$$
U\in e^+d^aC^{a-n/p}(\comega)+
C^{2a+0}({\Bbb R}^n)\cap  H_p^{2a}({\Bbb R}^n ),\tag2.18
$$ 
if $2a-n/p\ne 1$.
 If  $2a-n/p$ 
equals $1$, we need to add the space $\dot C^{1-0}(\comega)$.  

\endproclaim

\demo{Proof} The intersection with $\ol
H_p^{2a}(\complement\comega)$ serves as a bound at $\infty $. We extend $g$ to a function  $G\in
C^{2a+0}({\Bbb R}^n)$, then $G\in C^{2a+0}({\Bbb R}^n)\cap  H_p^{2a}({\Bbb
R}^n )$ (since $C^{t+0}\subset H^t_p$ over bounded sets). 
Theorem 2.3 now gives the existence of a solution $U=u+G$, where
$u\in H_p^{a(2a)}(\comega)$. By \cite{G15a} Cor.\ 5.5, cf.\ (2.9) above, this is
contained in   $ d^aC^{a-n/p}(\comega)$ when $2a-n/p\ne 1$ ($a-1/p$ and
$a-n/p$ are already noninteger). If $2a-p/n=1$, then we have to add the space $\dot
C^{1-0}(\comega)$, due 
to the embedding $\dot
H_p^{1+n/p}(\comega)\subset \dot C^{1-0}(\comega)$. 
\qed
\enddemo

Results for problems with $f\in L_\infty (\Omega )$ or H\"older-spaces
 were obtained in \cite{G15a} by letting $p\to\infty $; here we shall obtain
sharper results by applying the general method to the $C^s_*$-scale. Repeating the proof of Theorem 2.3 in this
scale, we find:

\proclaim{Theorem 2.5} Let $s>a-1$, and let $f\in\ol
C_*^{s-2a}(\Omega )$ and  $g\in \ol C^s_*(\complement\comega)$ be
given. Then the problem {\rm (2.15)} has the unique solution $U=u+G\in
C_*^{a(s)}(\comega)+C_*^s({\Bbb R}^n)$, where
$G\in C_*^s({\Bbb R}^n)$ is an extension of $g$ and 
$$
u=R(f-r^+P_aG)\in C_*^{a(s)}(\comega);\tag2.19
$$
here $R$ is the solution operator for {\rm (2.1)}.
\endproclaim

Let us spell this out in more detail for $s=2a$ and $s=2a+t$, in terms
of ordinary H\"older spaces. In the first statement, we take $g$ to be
compactly supported in $\complement\Omega $; in the next statements,
a very general term supported away from $\comega$ is added (it can in
particular lie in $C^{2a+t}_*$). Recall from (1.4) that $L_\infty
\subset C^0_*$.

\proclaim{Corollary 2.6}
$1^\circ$ For $f\in L_\infty (\Omega )$, $g\in
C_{\operatorname{comp}}^{2a}(\complement\Omega )$, the
solution of {\rm (2.15)} according to Theorem {\rm 2.5} satisfies 
$$
U\in e^+d^aC^{a}(\comega)\cap C^{2a}(\Omega )+
C_{\operatorname{comp}}^{2a}({\Bbb R}^n),\tag2.20
$$ 
with $2a$ replaced
by $1-0$ if $a=\frac12$. 

$2^\circ$ Let $X$ be any of the function spaces $F^\sigma _{p,q}({\Bbb
R}^n)$ or $B^\sigma _{p,q}({\Bbb R}^n)$, and denote by $X_{\operatorname{ext}}$
the subset of elements with support disjoint from $\comega$. For $f\in L_\infty (\Omega )$, $g\in
C_{\operatorname{comp}}^{2a}(\complement\Omega )+ X_{\operatorname{ext}}$, there
exists a solution $U$ of {\rm (2.15)} satisfying
$$
U\in e^+d^aC^{a}(\comega)\cap C^{2a}(\Omega )+
C_{\operatorname{comp}}^{2a}({\Bbb R}^n)+ X_{\operatorname{ext}},\tag2.21
$$
with $2a$ replaced by $1-0$ if $a=\frac12$. 

$3^\circ$ For $f\in C^t(\comega)$, $g\in
C^{2a+t}_{\operatorname{comp}}(\complement\Omega
)+X_{\operatorname{ext}}$, $t>0$, the solution according to $2^\circ$ satisfies
$$
U\in e^+d^aC^{a+t}(\comega)\cap C^{2a+t}(\Omega )+
C_{\operatorname{comp}}^{2a+t}({\Bbb R}^n)+ X_{\operatorname{ext}},\tag2.22$$
with $a+t$ resp.\ $2a+t$ replaced by $a+t-0$ resp.\ $2a+t-0$ when
they hit an integer.

\endproclaim

\demo{Proof}
$1^\circ$.
That $g\in C_{\operatorname{comp}}^{2a}(\complement\Omega )$ means
that $g$ is in $C^{2a}$ over the closed set $\complement\Omega $ and
vanishes outside a large ball; it extends to a function $G\in
C_{\operatorname{comp}}^{2a}({\Bbb R}^n )$. Since
$C^{2a}_{\operatorname{comp}}({\Bbb R}^n)\subset
C^{2a}_{\operatorname{comp},* }({\Bbb R}^n)$, the construction in
Theorem 2.5 gives a solution $U=u+G$, where $u$ is as in (2.12).

$2^\circ$. The function spaces are as described e.g.\ in \cite{J96},
 with $\sigma \in{\Bbb R}$, $0<p,q\le \infty $ ($p<\infty $ in the
 $F$-case), and $\psi $do's are well-defined in these spaces. We
 write $g=g_1+g_2$, where   $g_1\in
C_{\operatorname{comp}}^{2a}(\complement\Omega )$ and $g_2\in
X_{\operatorname{ext}}$.
The problem (2.15) with $g$ replaced by $g_1$ has a solution $u_1+G_1$
as under $1^\circ$. 
For the problem (2.15) with $f$ replaced by 0 and $g$
replaced by $g_2$ we take $G_2=g_2$. Then $P_aG_2$ is $C^\infty $ on a
neighborhood of $\comega$ (by the pseudolocal property of
pseudodifferential operators, cf.\ e.g.\ \cite{G09}, p.\ 177), so the reduced problem has a solution
$u_2\in \E_a(\comega)$, and the given problem then has the solution
$u_2+g_2$.

The sum of the solutions $u_1+G_1+u_2+g_2$ solves (2.15) and lies in
the asserted space.

$3^\circ$ is shown in a similar way, using (2.13).
\qed

\enddemo

\example{Remark 2.7} Note that according to the corollary, the effect
on the solution over
$\comega$, of an exterior
contribution to $g$ supported at a distance from $\comega$, is only a term in
$\E_a(\comega )$.
\endexample

\subhead 2.2 A nonhomogeneous Dirichlet problem \endsubhead

For the nonhomogeneous Dirichlet problem (2.2), the crucial
observation that leads to its solvability is that we can identify 
$
\E_{a-1}(\comega)/\E_{a}(\comega)$ with $ C^\infty (\Sigma  )$ by use
of the mapping
$$
\gamma _{a-1,0}\colon u\mapsto {\Gamma (a)}(d(x)^{1-a}u)|_{\Sigma
}\equiv \Gamma (a)\gamma _0(d^{1-a}u).\tag2.23
$$
(The gamma-function is included for consistency in calculations of 
 Fourier transformations and Taylor expansions.)
Namely, 
using normal and
tangential coordinates $x=y'+y_n\vec n(y')$ on a tubular neighborhood
$U_\delta =\{y'+y_n\vec n(y')\mid y'\in \Sigma , |y_n|<\delta \}$
of $\Sigma $ (where $\vec n(y')$ denotes the interior normal at $y'$), we have for $v\in C^\infty (\comega)$ that  
$$
v(x)=v(y'+y_n\vec n)=v_0(y')+y_nw(x)\text{ on }U_\delta \cap\comega ,
$$
where $v_0\in C^\infty (\Sigma )$ is the restriction of $v$ to $\Sigma $ (also denoted
$\gamma _0v$), and $w$ is $C^\infty $ on $U_\delta \cap \comega$. Now
when $u\in \E_{a-1}(\comega)$ is written as $u=e^+\Gamma
(a)^{-1}d(x)^{a-1}v$ with $v\in C^\infty (\comega)$,
$d(x)$ taken as $y_n$ on $U_\delta $, then 
$$
u(x)=\Gamma (a)^{-1}d(x)^{a-1}v_0(y')+ \Gamma
(a)^{-1}d(x)^aw(x)\text{ on }U_\delta \cap \comega, 
\tag2.24
$$
where $\Gamma (a)^{-1}d(x)^aw
$ is as a function in $ \E_a(\comega)$. 
Here $v_0$ is determined uniquely from $v$ and hence $\gamma _{a-1,0}u$
is determined uniquely from $u$, and the
null-space of the mapping $u\mapsto \gamma _{a-1,0}u$ is
$\E_{a}(\comega )$. See also Section 5 of \cite{G15a}; there it is
moreover shown  that this
mapping,
$$
\gamma _{a-1,0}\colon \E_{a-1}(\comega)\to C^\infty (\Sigma )\text{
with null-space }\E_a(\comega ),
$$
extends to a continuous surjective mapping
$$
\gamma _{a-1,0}\colon H_p^{(a-1)(s)}(\comega)\to B_p^{s-a+1/p'} (\Sigma )\text{
with null-space }H_p^{a(s)}(\comega ),\text{ for }s>a-1/p'. 
\tag2.25
$$

Now since we have the bijectiveness of $r^+P_a$ in (2.6), we can
simply adjoin the mapping (2.25) and conclude the bijectiveness of 
$$
\pmatrix 
r^+P_a\\ \quad\\\gamma _{a-1,0}\endpmatrix\colon
H_p^{(a-1)(s)}(\comega)\simto 
\matrix \ol
H_p^{s-2a}(\Omega )\\ \times \\ B_p^{s-a+1/p'} (\Sigma )\endmatrix.\tag2.26
$$
This gives the unique solvability of the problem (2.2) in these
spaces. There is an inverse $$
\pmatrix R& K\endpmatrix=\pmatrix r^+P_a\\ \gamma
_{a-1,0}\endpmatrix^{-1},$$ 
where $R$ is
the inverse of (2.6) as
introduced above, and $K$ is a mapping going from $\Sigma $ to $\comega
$. (Further details in \cite{G15a} Section 6.) 

In $C^s_*$-spaces, we likewise have an extension of
the mapping $\gamma _{a-1,0}$:
$$
\gamma _{a-1,0}\colon C_*^{(a-1)(s)}(\comega)\to C_*^{s-a+1} (\Sigma )\text{
with null-space }C_*^{a(s)}(\comega ),\text{ for }s>a-1. 
\tag2.27
$$
Then the result is as follows (as a
special case of Theorem 3.2 below), with conclusions in H\"older
spaces:

\proclaim{Theorem 2.8} Let $s>a-1$.  
The mapping $\{r^+P_a,\gamma _{a-1,0}\}$ defines a bijection
$$
\{r^+P_a,\gamma _{a-1,0}\}\colon C_*^{(a-1)(s)}(\comega)\to \ol
C_*^{s-2a}(\Omega )\times C_*^{s-a+1}(\Sigma ).\tag2.28
$$

In particular, for any $f\in L_\infty (\Omega )$, $\varphi \in C^{a+1}(\Sigma )$, there exists a
unique solution $u$ of {\rm (2.2)} in $C^{(a-1)(2a)}_*(\comega)$; it satisfies
$$
u\in \cases  e^+d(x)^{a-1}C^{a+1}(\comega)+\dot C^{2a}(\comega),\text{ when }a\ne
\tfrac12,\\
 e^+d(x)^{-\frac12}C^{\frac32}(\comega)+\dot
C^{1-0}(\comega),
\text{ when }a=
\tfrac12.\endcases
 \tag2.29
$$
 For $f\in C^t(\comega)$, $\varphi \in C^{a+1+t}(\Sigma ) $, $t>0$, the solution satisfies
$$
u\in \cases e^+d(x)^{a-1}C^{a+1+t}(\comega)+\dot C^{2a+t}(\comega),\text{ when
}a+t\text{ and }2a+t\notin  {\Bbb N},\\
e^+d(x)^{a-1}C^{a+1+t-0}(\comega)+\dot C^{2a+t-0}(\comega),\text{ when
}a+t\in  {\Bbb N},\\
e^+d(x)^{a-1}C^{a+1+t}(\comega)+\dot C^{2a+t-0}(\comega),\text{ when
 }2a+t\in  {\Bbb N}.
\endcases \tag2.30
$$
\endproclaim

\demo{Proof} The bijectiveness holds in view of the bijectiveness in
Theorem 2.2, and (2.27). The implications (2.29) and (2.30) follow from (2.10) with $\mu
=a-1$, together with the embedding properties recalled in Section 1.
Note that since $a+1>2a$, there is no need to mention an intersection
with $C^{2a(+t)}(\Omega )$.
\qed
\enddemo

This gives  sharpening of Th.\ 7.4 in \cite{G15a}.
We moreover recall that as shown in \cite{G15a} Th.\ 7.1, for functions
$u\in H_p^{(a-1)(s)}(\overline\Omega)$ for some $s,p$ with $s>a-1/p'$,
$$
f\in C^\infty (\comega),\; \varphi \in  C^{\infty }(\Sigma )\iff
u\in \E_{a-1}(\comega ).\tag2.31
$$

Also for the nonhomogeneous Dirichlet problem, there exist
formulations where the support condition on $u$ is replaced by a
prescription of its value on $\complement\Omega $. Abatangelo
\cite{A14} considers problems of the type
$$
\cases r^+P_aU& =f\text{ on }\Omega ,\\
U&=g\text{ on }\complement \Omega ,\\
\gamma _{a-1,0}U &=\varphi \text{ on }\Sigma .
\endcases \tag2.32
$$
(The boundary condition in \cite{A14} takes the form of the third line when $\Omega
$ is a ball, but is described in a more general way for other domains.)

For
(2.32), let  $f,g,\varphi $ be given with
$$
\{f,g,\varphi\}\in  \ol H_p^{s-2a}(\Omega )\times \ol
H^s_p(\complement\comega)
\times B_p^{s-a+1/p'} (\Sigma ),\text{ with }s>a-1/p'.\tag 2.33
$$
Then we search for a solution $U$ in a Sobolev space over ${\Bbb
R}^n$ that allows defining $\gamma _{a-1,0}U$.

We want to take as $G$
an extension of $g$ to $H_p^s({\Bbb R}^n)$. If
$s>n/p$, such that $H_p^s({\Bbb R}^{n})\subset C^0({\Bbb R}^n)$, we
have that
$\gamma _{a-1,0}\colon G\mapsto \Gamma (a)\gamma _0(d(x)^{1-a}G)$ is
well-defined and gives 0 for  $G\in H_p^s({\Bbb R}^n)$ (since
$a<1$). If $s<1/p$, we can take $G$ as the extension by 0 on $\Omega $
(since $\ol H_p^s(\complement\comega )$ identifies with $\dot H_p^s(\complement\Omega)$
when $-1/p'<s<1/p$). If $1/p\le s\le n/p$, we can also use the extension by 0 and note
that the boundary value from $\Omega $ is zero, but $G$ is only in $H_p^{1/p-0}({\Bbb R}^n)$. Now $U_1=U-G$ must satisfy  
$$
\cases r^+P_aU_1&=f-r^+P_aG\text{ in }\Omega ,\\
\supp U_1&\subset\comega,\\
\gamma _{a-1,0}U_1&=\varphi .\endcases \tag 2.34
$$
We continue the analysis for $s\notin [1/p,n/p]$; if $s$ is given $>0$, this can be
achieved by taking $p$ sufficiently large.

Since $P_aG\in H_{p,\operatorname{loc}}^{s-2a}({\Bbb R}^n)$,  $f-r^+P_aG\in \ol
H_p^{s-2a}(\Omega )$. Hereby we have reduced the problem to the form (2.3), where
we have the solution operator $\pmatrix R&K\endpmatrix$, see (2.26)ff.
This implies that (2.32) has the solution
$$
U=R(f-r^+P_aG)+K\varphi +G\in H_p^{a(s)}(\comega)+H_p^{(a-1)(s)}(\comega)
+ H_p^s({\Bbb R}^n).\tag2.35
$$
It is unique, since zero data give a zero solution (as we know from
(2.15) in the case $\varphi =0$). Recall that
$H_p^{a(s)}(\comega)\subset H_p^{(a-1)(s)}(\comega)$.

This shows the first part of the following theorem.

\proclaim{Theorem 2.9} $1^\circ$ Let $s>a-1/p'$ (if $s>0$ assume moreover
that $s\notin [1/p,n/p]$), and let $f, g,\varphi $ 
be
given as in {\rm (2.33)}. Let $G\in H_p^s({\Bbb R}^n)$ be an extension
of $g$ (by zero if $s<1/p$).

The problem {\rm (2.32)} has the unique
solution {\rm (2.35)} in $H_p^{(a-1)(s)}(\comega)+H^s_p({\Bbb R}^n)$.

$2^\circ$ Let $s>a-1$, $s\ne 0$, and let $f, g,\varphi $ 
be given with
$$
\{f,g,\varphi\}\in  \ol C_*^{s-2a}(\Omega )\times \ol
C_*^s(\complement\comega)
\times C_*^{s-a+1} (\Sigma ).\tag 2.36
$$
Let $G\in C_*^s({\Bbb R}^n)$ be an extension
of $g$ (by zero if $s<0$).

The problem {\rm (2.32)} has the unique solution 
$$
U=R(f-r^+P_aG)+K\varphi +G\in C_*^{(a-1)(s)}(\comega)
+ C_*^s({\Bbb R}^n).\tag2.37
$$
\endproclaim

\demo{Proof} $1^\circ$ was shown above, and $2^\circ$ is shown in an
analogous way:

For $s>0$, the extension $G$ has a boundary value $\gamma _{a-1,0}G= \Gamma (a)\gamma _0(d^{1-a}G)=0$
since $G$ is continuous and $1-a>0$, and for $s<0$ the boundary value from $\Omega
$ is 0, since $G$ is extended by zero (using that there is an
identification between $\ol C^s_*(\complement\comega)$
and $\dot C^s_*(\complement\Omega )$ when $-1<s<0$).
We then apply Theorem 2.8 to $u=U-G$.\qed
\enddemo

This reduction allows a study of higher regularity of the
solutions. The treatment in \cite{A14} seems primarily directed
towards solutions for not very smooth data. The boundary of
$\Omega $ is only assumed $C^{1,1}$ there.

\example{Remark 2.10}
When $s>a+n/p$, we note that since $H_p^{a(s)}(\comega)\subset
e^+d(x)^aC^0(\comega)\subset C^0({\Bbb R}^n)$ (cf.\ (2.9) or \cite{G15a} Cor.\ 5.5), the solution
(2.35) is
the sum of a continuous function and the term $K\varphi \in
H_p^{(a-1)(s)}(\comega)$ that stems solely from the boundary value
$\varphi $.
To further describe $K\varphi $, consider a localized situation,
where $\Omega $ is replaced by $\rnp$ and $d(x)$ is replaced by
$x_n$, and $P_a$ is carried over to a similar operator $P$ (of type and
factorization index $a$). As shown in the proof of 
 \cite{G15a} Th.\ 6.5, the solution $K\varphi $ (in a parametrix sense)
 of 
$$
r^+Pu=0\text{ in }\rnp,\quad \gamma _{a-1,0}u=\varphi  \text{ at }x_n=0,
$$
 is of the form $K\varphi =z+w$, where
$$
z=K_{a-1,0}\varphi=\Xi ^{1-a}_+e^+K_0\varphi 
=e^+c_{a-1}x_n^{a-1}K_0\varphi ,\quad w=-Rr^+Pz\in
H^{a(s)}(\crnp)\subset C^0({\Bbb R}^n);
$$
here $K_0$ is the standard Poisson operator sending $\varphi \in
B_p^{s-a+1/p'}({\Bbb R}^{n-1})$ into
$$
K_0\varphi =\F^{-1}_{\xi \to x}(\hat\varphi (\xi ')([\xi ']+i\xi _n)^{-1})=\F^{-1}_{\xi '\to
x'}(\hat\varphi (\xi ')e^{-[\xi ']x_n})\in \ol H_p^{s-a+1}(\rnp),
$$
with $\gamma _0K_0\varphi =\varphi $ (cf.\ also Cor. 5.3 and the proof
of Th.\ 5.4 in \cite{G15a}). Then
$$
z=e^+c_{a-1}x_n^{a-1}K_0\varphi  \in e^+x_n^{a-1}
\ol H_p^{s-a+1}(\rnp)\subset e^+x_n^{a-1}C^{s-a+1-n/p}(\crnp),
$$
with $K_0\varphi \ne 0$ at $\{x_n=0\}$
 when $\varphi \ne 0$. For higher $s$, the factor $K_0\varphi $ lies
 in  higher-order Sobolev and H\"older spaces,
but is always 
nontrivial at $\{x_n=0\}$ when $\varphi \ne 0$.

When this is carried back to the manifold situation, we have that $U$ is the
sum of a term in $C^0({\Bbb R}^n)$ and a term  $ e^+d(x)^{a-1}v$, $v\in
\ol H_p^{s-a+1}(\Omega )$, where $v$ is nonzero at $\partial\Omega $ when $\varphi \ne0$.
Since $a<1$, this term blows up at the boundary.

Hence the solutions are ``large'' at the boundary in this precise
 sense, consisting of a continuous function plus a term containing the
 factor $d(x)^{a-1}$ nontrivially. Cf.\ also (2.31).

 It is a theme of \cite{A14} that there exist ``large'' solutions of
 the nonhomogeneous Dirichlet problem; we here see that this is {\it not
 an exception but a rule} of the setup, provided naturally by the
 part of the solution mapping going from $\Sigma $ to $\comega$. 
\endexample

Theorem 2.9 $1^\circ$ gives the following result in H\"older spaces when
$f\in L_p(\Omega )=\ol H^0_p(\Omega )$.

\proclaim{Corollary 2.11} 
Let $p>n/a$. For $f\in L_p(\Omega )$, $g\in
C^{2a+0}(\complement\Omega )\cap \ol H_p^{2a}(\complement\comega )$ and
$\varphi \in C^{a+1/p'+0}(\Sigma )$, the
solution $U$ of {\rm (2.32)} according to Theorem {\rm 2.8} satisfies
$$
U\in e^+d^{a-1}C^{a+1-n/p}(\comega)+\dot C^{2a-n/p}(\comega)+
C^{2a+0}({\Bbb R}^n)\cap  H_p^{2a}({\Bbb R}^n ),\tag2.38
$$ 
with $2a-n/p$ replaced by $1-0$ if $2a-n/p=1$.

\endproclaim

\demo{Proof} 
Note that $2a>n/p$. We extend $g$ as in Corollary 2.4 to a
function $G\in C^{2a+0}({\Bbb R}^n)\cap  H_p^{2a}({\Bbb R}^n ) $, and
note that $\varphi \in C^{a+1/p'+0}(\Sigma )\subset
B_p^{a+1/p'}(\Sigma )$. Theorem 2.9 $1^\circ$
shows that there is a (unique) solution $U=u+K\varphi +G$ with 
$$
u+K\varphi \in H_p^{(a-1)(2a)}(\comega)\subset
e^+d^{a-1}C^{a+1-n/p}(\comega)+\dot C^{2a-n/p}(\comega)
$$
(one may consult \cite{G15a} (7.12)), with the mentioned modification if $2a-n/p$
is integer.  
\qed
\enddemo

For $f\in L_\infty (\Omega )$ or $C^t(\comega)$, we get the sharpest
results by applying the statement for $C^s_*$-spaces:
 
\proclaim{Corollary 2.12}
$1^\circ$ For $f\in L_\infty (\Omega )$, $g\in
C_{\operatorname{comp}}^{2a}(\complement\Omega )$ and $\varphi \in C^{a+1}(\Sigma )$, the
solution of {\rm (2.32)} satisfies 
 $$
U\in e^+d^{a-1}C^{a+1}(\comega)+
C_{\operatorname{comp}}^{2a}({\Bbb R}^n),\tag2.39
$$ 
with $2a$ replaced
by $1-0$ if $a=\frac12$. 

$2^\circ$ Let $X$ be any of the function spaces $F^\sigma _{p,q}({\Bbb
R}^n)$ or $B^\sigma _{p,q}({\Bbb R}^n)$, and denote by $X_{\operatorname{ext}}$
the subset of elements with support disjoint from $\comega$. For $f\in L_\infty (\Omega )$, $g\in
C_{\operatorname{comp}}^{2a}(\complement\Omega )+
X_{\operatorname{ext}}$ and $\varphi \in C^{a+1}(\Sigma )$, there
exists a solution of {\rm (2.32)} satisfying 
$$
U\in e^+d^{a-1}C^{a+1}(\comega)+
C_{\operatorname{comp}}^{2a}({\Bbb R}^n)+ X_{\operatorname{ext}},\tag2.40
$$
with $2a$ replaced by $1-0$ if $a=\frac12$. 

$3^\circ$ For $f\in C^t(\comega)$, $g\in
C^{2a+t}_{\operatorname{comp}}(\complement\Omega
)+X_{\operatorname{ext}}$ and $\varphi \in C^{a+1+t}(\Sigma )$, the solution according to $2^\circ$ satisfies
$$
U\in e^+d^{a-1}C^{a+1+t}(\comega)+
C_{\operatorname{comp}}^{2a+t}({\Bbb R}^n)+ X_{\operatorname{ext}},$$
with $a+1+t$ resp.\ $2a+t$ replaced by $a+1+t-0$ resp.\ $2a+t-0$ when
they hit an integer.

\endproclaim

\demo{Proof} We apply Theorem 2.9 $2^\circ$ very much in the same way
as in Corollary 2.6; details can be omitted.
\qed
\enddemo

\subhead 2.3 A nonhomogeneous Neumann problem \endsubhead

The Neumann boundary value defined in connection with $(-\Delta )^a$ is 
$$
\gamma _{a-1,1}u=\Gamma (a+1)\gamma _0(\partial_n(d(x)^{1-a}u));\tag2.41
$$
it is proportional to the second coefficient in the Taylor expansion
of $d^{1-a}u$ in the normal variable at the boundary (like $\gamma _0w$ when $w$ is as in (2.24)).

We here have, by use of Theorem 4.3 below:

\proclaim{Theorem 2.13} The mapping $\{r^+P_a,\gamma _{a-1,1}\}$ defines a
Fredholm operator: 
$$
\{r^+P_a,\gamma _{a-1,1}\}\colon H_p^{(a-1)(s)}(\comega)\to \ol H_p^{s-2a}(\comega)\times B_p^{s-a-1/p}(\Sigma ),\tag 2.42
$$
for $s>a+1/p$.

\endproclaim
 \demo{Proof} The continuity of the mapping (2.42) follows from
 \cite{G15a} Th.\ 5.1 with $\mu =a-1$, $M=2$. The Fredholm property
 follows from Theorem 4.3 below in a special case, cf.\ (3.2), by piecing together a parametrix from
 the parametrix construction in local coordinates given there. We use
 that the parametrix exists since $P_a$ in local coordinates has
 principal symbol $|\xi |^{2a}$. \qed 
\enddemo

There is a similar  version in $C^s_*$-spaces, with
consequences for H\"older estimates:

\proclaim{Theorem 2.14} Let $s>a$.  
The mapping $\{r^+P_a,\gamma _{a-1,1}\}$ defines a Fredholm operator
$$
\{r^+P_a,\gamma _{a-1,1}\}\colon C_*^{(a-1)(s)}(\comega)\to \ol
C_*^{s-2a}(\Omega )\times C_*^{s-a}(\Sigma ).\tag2.43
$$

In particular, for  $\{f,\psi \}\in L_\infty (\Omega )\times
C^{a}(\Sigma )$ subject to a certain finite set of linear constraints there exists a
solution $u$ of {\rm (2.3)} in $C^{(a-1)(2a)}_*(\comega)$; it is unique modulo a finite dimensional linear subspace
$\Cal N\subset \E_{a-1}(\comega)$ and satisfies
$$
u\in \cases  e^+d(x)^{a-1}C^{a+1}(\comega)+\dot C^{2a}(\comega),\text{ when }a\ne
\tfrac12,\\
 e^+d(x)^{-\frac12}C^{\frac32}(\comega)+\dot
C^{1-0}(\comega),
\text{ when }a=
\tfrac12.\endcases
 \tag2.44
$$
 For $f\in C^t(\comega)$, $\psi  \in C^{a+t}(\Sigma ) $, $t>0$, the solution satisfies
$$
u\in \cases e^+d(x)^{a-1}C^{a+1+t}(\comega)+\dot C^{2a+t}(\comega),\text{ when
}a+t\text{ and }2a+t\notin  {\Bbb N},\\
e^+d(x)^{a-1}C^{a+1+t-0}(\comega)+\dot C^{2a+t-0}(\comega),\text{ when
}a+t\in  {\Bbb N},\\
e^+d(x)^{a-1}C^{a+1+t}(\comega)+\dot C^{2a+t-0}(\comega),\text{ when
 }2a+t\in  {\Bbb N}.
\endcases \tag2.45
$$
\endproclaim

\demo{Proof} The first statement is the analogue of Theorem 2.13,
now derived from Theorem 4.3 for $p=q=\infty $.
In the next, detailed statements we formulate the Fredholm
property explicitly, using also Theorem 3.5 on the smoothness of the
kernel. Here the inclusions (2.44) and (2.45) follow from the
description (2.10)  of $C^{(a-1)(s)}_*(\comega)$ as in the
proof of Theorem 2.8.\qed 
\enddemo

Also in the Neumann case, one can formulate versions of the theorems with $u$
prescribed on ${\Bbb R}^n\setminus \Omega $, and show their
equivalence with the set-up for $u$ supported in $\comega$; we think
this is sufficiently exemplified by the treatment of the Dirichlet
condition above, that we
can leave details to the interested reader.

\head 3. Boundary problems in general spaces
\endhead

One of the conclusions in \cite{G15a} of the study of the $\psi $do $P$ of order $m\in{\Bbb C}$,
with factorization index and
type $\mu _0\in{\Bbb C}$, was that it could be linked, by the help of
the special order-reducing operators $\Lambda _\pm^{(\mu )}$, to an operator
$$
Q=\Lambda _-^{(\mu _0-m)}P\Lambda _+^{(-\mu _0)}\tag3.1
$$
of order 0 and with factorization index and type 0, which could be
treated by the help of the
calculus of Boutet de Monvel on $H^s_p$-spaces, as accounted
for in \cite{G90}. Results for $P$ and its boundary value
problems could then be deduced from those for $Q$ in the case of a
homogeneous boundary condition. With a natural definition of boundary
operators $\gamma _{\mu ,k}$, also nonhomogeneous boundary conditions
could be treated. In particular, we found the structure of
parametrices of $r^+P$, with homogeneous or nonhomogeneous Dirichlet-type
conditions, as compositions of operators belonging to the Boutet de
Monvel calculus with the special order-reducing operators, see Theorems
4.4, 6.1 and 6.5 of \cite{G15a}.

The results of \cite{G90} have been extended to the much more general
families of  spaces 
$F^s_{p,q}$ (Triebel-Lizorkin spaces) and $B^s_{p,q}$ (Besov spaces)
by
Johnsen in \cite{J96}. He shows that elliptic systems on a compact manifold with a smooth boundary, belonging to the Boutet de
Monvel calculus, have Fredholm solvability also in these more general
spaces, with $C^\infty $ kernels and range complements (cokernels) independent of $s,p,q$.
Here $0<p,q\le\infty $ is allowed for the $B^s_{p,q}$-spaces, and the
same goes for the $F^s_{p,q}$-spaces, except that $p$ is taken $< \infty$  
(to avoid long explanations of exceptional cases). The parameter
$s$ is taken $>s_0$, for a suitable $s_0$ depending on $p$ and the order and class
of the involved operators. We refer to \cite{J96} (or to Triebel's
books) for  detailed descriptions
of the spaces, just recalling that for $1<p<\infty $,
$$
\aligned
F^s_{2,2}&=B^s_{2,2}=H^s_2,\text{ $L_2$-Sobolev spaces,}\\
F^s_{p,2}&=H^s_p,\text{ Bessel-potential spaces,}\\
B^s_{p,p}&=B^s_p,\text{ Besov spaces.} 
\endaligned\tag3.2
$$
Here the Bessel-potential spaces $H^s_p$ are also called $W^s_p$ (or
$W^{s,p}$) for
$s\in{\Bbb N}_0$, and the Besov spaces $B^s_p$ are also called $W^s_p$ (or
$W^{s,p}$) for
$s\in\rp\setminus {\Bbb N}$, under the common name
Sobolev-Slobodetskii spaces. Recall moreover that
$F^s_{p,p}=B^s_{p,p}$ for $0<p<\infty $ (also denoted $B^s_p$).

We return to the general situation of $\comega $ smoothly embedded in a
Riemannian manifold $\Omega _1$, with $\crnp\subset{\Bbb R}^n$ used in
localizations.  H\"ormander's notation $\dot F, \ol F$
and $\dot B, \ol B$ will be used for the general scales, in the same
way as for $H^s_p$, cf.\ (1.2)ff. 

In the present paper, we shall in particular be interested in the case of the scale of
spaces $B^s_{\infty ,\infty }=C_*^s$ (see the text around (1.3)), which gives a shortcut to sharp results on
solvability in H\"older spaces.

{\it Since we are mostly interested in results for large $p$, we shall
assume $p\ge 1$}, which simplifies the quotations from \cite{J96}, namely,
the condition $s>\max\{1/p-1,n/p-n\}$ simplifies to $s>1/p-1$, since
$1/p-1\ge n/p-n$ when $p\ge 1$. (In situations
where $p<1$ would be needed, e.g.\ in bootstrap regularity arguments,
one can supply the presentation here with the appropriate results from \cite{J96}.)  
The usual notation $1/p'=1-1/p$ is understood as $0$ resp.\ $1 $ when $p=1$
resp.\ $\infty $. We assume $p\le \infty $ in  $B$-cases,
$p<\infty $ in $F$-cases, and take $0< q\le \infty $.

The scales $F^s_{p,q}$ and $B^s_{p,q}$ have analogous roles in
definitions over $\comega$, but the 
trace mappings on them are slightly different: When $s>1/p$, 
$$
\gamma _0\colon \ol F^s_{p,q}(\Omega )\to B^{s-1/p}_{p,p}(\partial\Omega
),\quad \gamma _0\colon \ol B^s_{p,q}(\Omega )\to B^{s-1/p}_{p,q}(\partial\Omega ),\tag3.3
$$
continuously and surjectively. (One could also write $F^s_{p,p}$
instead of $B^s_{p,p}$; in \cite{J96}, both indications occur.) 

To reduce repetitive formulations, we shall introduce the common
notation:
$$
X^{s}_{p,q} \text{ stands for either }F^{s}_{p,q}\text{ or }B^{s}_{p,q},\text{
at convenience, }\tag3.4
$$
with the same choice in each place if the notation appears several times in the same
calculation. Formulas involving boundary operators will be given
explicitly in the
two different cases resulting from (3.3).

In addition to the mapping and Fredholm properties established for
Boutet de Monvel systems in \cite{J96}, we need the following generalizations of
(1.6) (as in \cite{G15a} (1.11)--(1.20)): 
$$
\aligned
\Xi ^\mu _+\text{ and }\Lambda^\mu _+&\colon \dot X_{p,q}^s(\crnp)\simto
\dot X_{p,q}^{s-\operatorname{Re}\mu }(\crnp),\text{ with
inverse }\Xi _+^{-\mu }\text{ resp.\ }\Lambda
^{-\mu }_+,\\
\Xi ^\mu _{-,+}\text{ and }\Lambda ^\mu _{-,+}&\colon \ol X_{p,q}^s(\rnp)\simto
\ol X_{p,q}^{s-\operatorname{Re}\mu }(\rnp),\text{ with
inverse }\Xi ^{-\mu }_{-,+}\text{ resp.\ }\Lambda ^{-\mu }_{-,+},\\
\Lambda^{(\mu )}_+&\colon \dot X_{p,q}^s(\comega )\simto
\dot X_{p,q}^{s-\operatorname{Re}\mu }(\comega ),\\
\Lambda ^{(\mu )}_{-,+}&\colon \ol X_{p,q}^s(\Omega )\simto
\ol X_{p,q}^{s-\operatorname{Re}\mu }(\Omega ),
\endaligned
\tag3.5
$$ 
valid for all $s\in{\Bbb R}$. The cases with integer $\mu $ are
covered by \cite{J96} as a direct extension of the presentation in
\cite{G90}, the cases of more general $\mu $ likewise extend, since
the support preserving properties extend. 

We can then define (analogously to the definitions and observations in \cite{G15a},
Sect.\ 1.2, 1.3):

\proclaim{Definition 3.1} Let $s>\operatorname{Re}\mu -1/p'$. 

$1^\circ$ A distribution $u$ on ${\Bbb R}^n$ is in $X_{p,q}^{\mu
(s)}(\crnp)$ if and only if $\Xi _+^\mu u\in \dot X_{p,q}^{-1/p'+0}(\crnp)$
and $r^+\Xi _+^\mu u\in \ol X_{p,q}^{s-\operatorname{Re}\mu }(\rnp)$. In
fact, $r^+\Xi _+^\mu $ maps $X_{p,q}^{\mu (s)}(\crnp)$ bijectively onto
$\ol X_{p,q}^{s-\operatorname{Re}\mu }(\rnp)$ with inverse $\Xi _+^{-\mu
}e^+$, and
$$
X_{p,q}^{\mu (s)}(\crnp)=\Xi _+^{-\mu }e^+\ol X_{p,q}^{s-\operatorname{Re}\mu }(\rnp),\tag3.6
$$
with the inherited norm. Here $\Lambda ^{-\mu }_+$ can equivalently be
used.

$2^\circ$ A distribution $u$ on $\Omega _1$ is in $X_{p,q}^{\mu
(s)}(\comega)$ if and only if $\Lambda _+^{(\mu )}u\in \dot X_{p,q}^{-1/p'+0}(\comega)$
and $r^+\Lambda  _+^{(\mu )}u\in \ol X_{p,q}^{s-\operatorname{Re}\mu }(\Omega )$. In
fact, $r^+\Lambda  _+^{(\mu )}$ maps $X_{p,q}^{\mu (s)}(\comega)$ bijectively onto
$\ol X_{p,q}^{s-\operatorname{Re}\mu }(\Omega )$ with inverse $\Lambda  _+^{(-\mu
)}e^+$, and
$$
X_{p,q}^{\mu (s)}(\comega)=\Lambda _+^{(-\mu )}e^+ \ol X_{p,q}^{s-\operatorname{Re}\mu }(\Omega),\tag3.7
$$
with the inherited norm.
\endproclaim

The distributions in  $X_{p,q}^{\mu (s)}(\crnp)$ resp.\ $X_{p,q}^{\mu
(s)}(\comega)$ are locally in $X_{p,q}^s$ over $\rnp$ resp.\ $\Omega $, by
interior regularity.

By use of the mapping properties of the standard trace operators
$\gamma _j$ described in \cite{J96}, and use of (3.5) above, the trace
operators $\varrho _{\mu ,M}$ introduced in \cite{G15a}, Sect.\ 5, extend to the general spaces:
$$
\varrho _{\mu ,M}=\{\gamma _{\mu ,0}, \gamma _{\mu ,1},\dots,\gamma _{\mu ,M-1}\}\colon \cases F_{p,q}^{\mu (s)}(\comega)\to  \prod _{0\le
j<M}B_{p,p}^{s-\operatorname{Re}\mu -j-1/p}(\partial\Omega ),\\
B_{p,q}^{\mu (s)}(\comega)\to  \prod _{0\le
j<M}B_{p,q}^{s-\operatorname{Re}\mu -j-1/p}(\partial\Omega ),
\endcases\tag3.8
$$
for $s>\operatorname{Re}\mu +M-1/p'$;
they are surjective with kernels $F_{p,q}^{(\mu +M)(s)}(\comega)$ and $B_{p,q}^{(\mu +M)(s)}(\comega)$.

We can now formulate some important results from \cite{G15a} in these
scales of spaces. Recall that when $P$ is of type $\mu $, it is also of
type $\mu '$ for $\mu -\mu '\in{\Bbb Z}$.

\proclaim{Theorem 3.2}

$1^\circ$ Let the $\psi $do $P$ on $\Omega _1$ be of order $m\in{\Bbb
C}$ and of type $\mu \in{\Bbb C}$ relative to the boundary of the
smooth compact subset  $\comega
\subset \Omega _1$. Then when $s>\operatorname{Re}\mu -1/p'
$, $r^+P$ maps $X_{p,q}^{\mu (s)}(\comega)$ continuously into $\ol
X_{p,q}^{s-\operatorname{Re}m}(\Omega)$.

$2^\circ$ Assume in addition that $P$ is elliptic 
 and has factorization index $\mu _0$, where $\mu - \mu _0\in{\Bbb Z}$. Let
$s>\operatorname{Re}\mu _0-1/p'$. If $u\in \dot X_{p,q}^{\sigma }(\comega)$ for some $\sigma >\operatorname{Re}\mu _0-1/p'$ and
$r^+Pu\in \ol
X_{p,q}^{s-\operatorname{Re}m}(\Omega)$, then $u\in X_{p,q}^{\mu _0(s)}(\comega)$. The mapping
$r^+P$ defines a Fredholm operator
$$
r^+P\colon X_{p,q}^{\mu _0(s)}(\comega)\to \ol
X_{p,q}^{s-\operatorname{Re}m}(\Omega).\tag3.9
$$
Moreover, $\{r^+P,\gamma _{\mu _0-1,0}\}$ defines a Fredholm operator
$$
\{r^+P,\gamma _{\mu _0-1,0}\}\colon \cases F_{p,q}^{(\mu _0-1)(s)}(\comega)\to \ol
F_{p,q}^{s-\operatorname{Re}m}(\Omega )\times
B_{p,p}^{s-\operatorname{Re}\mu _0+1-1/p}(\partial\Omega ),\\
B_{p,q}^{(\mu _0-1)(s)}(\comega)\to \ol
B_{p,q}^{s-\operatorname{Re}m}(\Omega )\times
B_{p,q}^{s-\operatorname{Re}\mu _0+1-1/p}(\partial\Omega ).
\endcases\tag3.10
$$

$3^\circ$ Let $P$ be as in $2^\circ$, and let $\mu '=\mu _0-M$ for a
positive integer $M$. Then when $s>\operatorname{Re}\mu _0-1/p'$, $\{r^+P,\varrho _{\mu ',M}\}$ 
defines a Fredholm operator
$$
\{r^+P,\varrho _{\mu ',M}\}\colon 
\cases F_{p,q}^{\mu '(s)}(\comega)\to \ol
F_{p,q}^{s-\operatorname{Re}m}(\Omega )\times \prod _{0\le
j<M}B_{p,p}^{s-\operatorname{Re}\mu '-j-1/p}(\partial\Omega ),\\
B_{p,q}^{\mu '(s)}(\comega)\to \ol
B_{p,q}^{s-\operatorname{Re}m}(\Omega )\times \prod _{0\le
j<M}B_{p,q}^{s-\operatorname{Re}\mu '-j-1/p}(\partial\Omega ).
\endcases\tag3.11
$$
    
\endproclaim 

\demo{Proof} $1^\circ$. This is the extension of \cite{G15a} Th.\ 4.2
to the general spaces. We recall that the proof consist of a reduction
of the study of $r^+P$  to the
consideration of $Q_+$ (with $Q$ as in (3.1) for $\mu =\mu _0$) of
type 0; this works well in the present spaces.

$2^\circ$--$3^\circ$. Here (3.9) is obtained by a generalization of \cite{G15a}
Th.\ 4.4 and its proof to the current
spaces. Now (3.11) is obtained as in \cite{G15a} Th.\ 6.1 by adjoining the mapping
(3.8) (with $\mu =\mu '$) to $r^+P$. Here (3.10) is the special case $M=1$, as in
\cite{G15a} Cor.\ 6.2.\qed

\enddemo

The parametrices $R$ and $\pmatrix R&K\endpmatrix$ described by formulas in \cite{G15a} Th.\ 4.4 and 6.5 also
work in the present spaces.

\example{Example 3.3} As an example, we have for the choice $X=B$,
$p=q=\infty $, i.e., 
$X^s_{p,q}=B^s_{\infty ,\infty }=C^s_*$, that Theorem 3.2 $2^\circ$ shows the following:

Let $P$ be elliptic of order $m$ and of type $\mu _0$, with
factorization index $\mu _0$, and let
$s>\operatorname{Re}\mu _0-1$. If $u\in \dot C_*^{\sigma }(\comega)$ for some $\sigma >\operatorname{Re}\mu _0-1$ and
$r^+Pu\in \ol
C_*^{s-\operatorname{Re}m}(\Omega)$, then $u\in C_*^{\mu _0(s)}(\comega)$. The mapping
$r^+P$ defines a Fredholm operator
$$
r^+P\colon C_*^{\mu _0(s)}(\comega)\to \ol
C_*^{s-\operatorname{Re}m}(\Omega).\tag3.12
$$
Moreover, $\{r^+P,\gamma _{\mu _0-1,0}\}$ defines a Fredholm operator
$$
\{r^+P,\gamma _{\mu _0-1,0}\}\colon  C_*^{(\mu _0-1)(s)}(\comega)\to \ol
C_*^{s-\operatorname{Re}m}(\Omega )\times
C_*^{s-\operatorname{Re}\mu _0+1}(\partial\Omega ).
\tag3.13
$$

\endexample

For $\operatorname{Re}\mu >-1/p'$, the spaces $X_{p,q}^{\mu (s)}(\crnp)$ 
and $X_{p,q}^{\mu (s)}(\comega)$ are further described by the
following generalization of \cite{G15a}, Th.\ 5.4:

\proclaim{Theorem 3.4}
One has for
$\operatorname{Re}\mu >-1$, $s>\operatorname{Re}\mu -1/p'$, with $M\in{\Bbb N}$:
$$
\aligned
&X_{p,q}^{\mu (s)}(\crnp)\cases =\dot X_{p,q}^s(\crnp)\text{ if
}s-\operatorname{Re}\mu \in \,]-1/p',1/p[\,,\\
\subset\dot X_{p,q}^{s-0}(\crnp)\text{ if }s-\operatorname{Re}\mu =1/p.\endcases
\\
X_{p,q}^{\mu (s)}(\crnp)&\subset e^+x_n^{\mu }\ol X_{p,q}^{s-\operatorname{Re}\mu
}(\rnp)+
\cases \dot X_{p,q}^s(\crnp)\text{ if }s-\operatorname{Re}\mu \in M+\,]-1/p',1/p[\, \\
\dot X_{p,q}^{s-0}(\crnp)\text{ if }s-\operatorname{Re}\mu = M+1/p.\endcases
\endaligned\tag 3.14
 $$

The inclusions {\rm (3.14)} also hold in the manifold situation, with
$\rnp$ replaced by $\Omega $ and $x_n$ replaced by $d(x)$.

\endproclaim

\demo{Proof} The first statement in (3.14) follows since
$e^+\ol X_{p,q}^t(\rnp)=\dot X_{p,q}^t(\crnp)$ for
$-1/p'<t<1/p$, cf.\ \cite{J96} (2.51)--(2.52).

For the second statement we use the representation of $u$ as in \cite{G15a}
(5.13)--(5.14), in the same way as in the proof of Th.\ 5.4 there. The
crucial fact is that the Poisson operator $K_0$ maps $\gamma _{\mu
,0}u\in B_{p,p}^{s-\operatorname{Re}\mu -1/p}({\Bbb R}^{n-1})$ resp.\  $B_{p,q}^{s-\operatorname{Re}\mu -1/p}({\Bbb R}^{n-1})$
into $\ol F_{p,q}^{s-\operatorname{Re}\mu }(\rnp)$ resp.\ $\ol B_{p,q}^{s-\operatorname{Re}\mu }(\rnp)$ 
(by
\cite{J96}), defining a term 
$$
v_0=e^+K_{\mu ,0}\gamma _{\mu ,0}u=c_\mu e^+x_n^\mu K_0\gamma _{\mu,0}u\in e^+x_n^\mu\ol X_{p,q}^{s-\operatorname{Re}\mu }(\rnp),
$$ with similar descriptions of terms $e^+K_{\mu ,j}\gamma _{\mu
,j}u$ for $j$ up to $M-1$, such that $u$ by subtraction of these terms
gives a term in $\dot X_{p,q}^s(\crnp)$ (with $s$ replaced
by $s-0$ if $s-\operatorname{Re}\mu -1/p$ hits an integer).
 \qed
\enddemo

Moreover, it is important to observe the following invariance property of kernels and
cokernels (typical in elliptic theory):

\proclaim{Theorem 3.5} For the Fredholm operators considered in
Theorem {\rm 3.2}, the kernel is a finite dimensional subspace $\Cal N$ of $\E_\mu
(\comega)$ independent of the choice  of $s,p,q$ and $F$ or $B$. 

There is a finite dimensional
range complement $\Cal M\subset C^\infty (\comega)$ for {\rm
(3.9)}, resp.\
$\Cal M_1\subset C^\infty (\comega)\times C^\infty (\partial\Omega )^M
$ for {\rm (3.10)--(3.11)}, that is
independent of  the choice  of $s,p,q$, $F$, $B$. 
\endproclaim

\demo{Proof} This follows from the similar statement for operators in the
Boutet de Monvel calculus in \cite{J96} Sect.\ 5.1, when we apply the
mappings $\Lambda ^{(\mu )} _\pm$ etc.\ in the reduction of the homogeneous
Dirichlet problem to a problem in the Boutet de Monvel calculus. \qed
\enddemo

\head 4. More general boundary conditions \endhead

In Theorem 3.2 we have obtained the Fredholm solvability of Dirichlet-type problems
defined by operators
$$
\{r^+P,\gamma _{\mu -1 ,0}\}\colon \cases F_{p,q}^{(\mu -1)(s)}(\comega)\to \ol
F_{p,q}^{s-\operatorname{Re}m}(\Omega )\times
B_{p,p}^{s-\operatorname{Re}\mu +1/p'}(\partial\Omega ),\\
B_{p,q}^{(\mu -1)(s)}(\comega)\to \ol
B_{p,q}^{s-\operatorname{Re}m}(\Omega )\times
B_{p,q}^{s-\operatorname{Re}\mu +1/p'}(\partial\Omega ),
\endcases\tag4.1
$$
for $s>\operatorname{Re}\mu -1/p'$, where $P$ is elliptic of order
$m$, is of
type $\mu $, and
has factorization index $\mu $ (called $\mu _0$ there). In Th.\ 6.5 of \cite{G15a} we constructed a parametrix
in local coordinates, which in the Besov-Triebel-Lizorkin scales maps as follows:
$$
\pmatrix R_D& K_D\endpmatrix \colon
\cases 
\ol
F_{p,q}^{s-\operatorname{Re}m}(\rnp)\times
B_{p,p}^{s-\operatorname{Re}\mu +1/p'}({\Bbb R}^{n-1} )\to
F_{p,q}^{(\mu -1)(s)}(\crnp),\\
\ol
B_{p,q}^{s-\operatorname{Re}m}(\rnp)\times
B_{p,q}^{s-\operatorname{Re}\mu +1/p'}({\Bbb R}^{n-1} )\to
B_{p,q}^{(\mu -1)(s)}(\crnp),
\endcases\tag4.2
$$
where $R_D=\Lambda ^{-\mu }_+e^+  \widetilde {Q_+}\Lambda ^{\mu
-m}_{-,+}$ and  $K_D=\Xi _+^{1-\mu }e^+K'$ or $\Lambda _+^{1-\mu }e^+K''
$, $\widetilde {Q_+}$ being a parametrix of $Q_+$ (where $Q$ is
recalled in (3.1)), and $K'$ and $K''$ being Poisson operators in the
Boutet de Monvel calculus of order 0.

\subhead 4.1 Boundary operators of type $\gamma _0B$ \endsubhead

We shall now describe a general way to let other 
boundary operators enter in lieu
of $\gamma _{\mu -1,0}$.  The point is to reduce the problem to a
problem in the Boutet de Monvel calculus (with $\psi $do's of type 0
and integer order). We can assume that the family of  auxiliary operators
$\Lambda _{\pm}^{(\varrho )}$ is chosen such that $(\Lambda _{\pm}^{(\varrho )})^{-1}=\Lambda _\pm^{(-\varrho )}$.
  
\proclaim{Theorem 4.1} Let $P$ be elliptic of order $m\in {\Bbb C}$ on
 $\Omega _1$,
having  type $\mu $ and 
factorization index $\mu $ with respect to the smooth compact subset
 $\comega$. Let $B$ be a $\psi $do of order $m_0+\mu $ and of type
$\mu $, with $m_0$ integer. Consider the mapping
$$
\{r^+P,\gamma _{0}r^+B\}\colon \cases F_{p,q}^{(\mu -1)(s)}(\comega)\to \ol
F_{p,q}^{s-\operatorname{Re}m}(\Omega )\times
B_{p,p}^{s-m_0-\operatorname{Re}\mu +1/p'}(\partial\Omega ),\\
B_{p,q}^{(\mu -1)(s)}(\comega)\to \ol
B_{p,q}^{s-\operatorname{Re}m}(\Omega )\times
B_{p,q}^{s-m_0-\operatorname{Re}\mu +1/p'}(\partial\Omega ),
\endcases
\tag4.3
$$
for $s>\operatorname{Re}\mu +\max\{m_0,0\}-1/p'$.

$1^\circ$ For $u\in X_{p,q}^{(\mu -1)(s)}(\crnp)$, the problem
$$
r^+Pu=f\text{ on }\Omega ,\quad \gamma _0r^+Bu=\psi\text{ on
}\partial\Omega  , \tag4.4
$$
can be reduced to an equivalent problem
$$
P'_+w=g\text{ on }\Omega ,\quad \gamma _0B'_+w=\psi\text{ on }\partial\Omega ,\tag4.5
$$
where $w=r^+\Lambda _+^{(\mu -1)}u\in \ol X_{p,q}^{s-\operatorname{Re}\mu
+1}(\Omega )$,  $g=\Lambda ^{(\mu -m)}_{-,+}f\in \ol X_{p,q}^{s-\operatorname{Re}\mu
}(\Omega )$, and where
$$
P'=\Lambda _-^{(\mu -m)}P\Lambda _{+}^{(1-\mu )},\quad
B'=B\Lambda _{+}^{(1-\mu )},\tag4.6$$
$\psi $do's of order $1$ resp.\ $m_0 +1$, and type $0$. 

$2^\circ$ The problem {\rm (4.4)} is Fredholm solvable for $s>\operatorname{Re}\mu +\max\{m_0,0\}-1/p'$, if and only if the
problem {\rm (4.5)} is Fredholm solvable, as a mapping
$$
\{P'_+,\gamma _{0}B'_+\}\colon \cases\ol F_{p,q}^{t +1}(\Omega)\to \ol
F_{p,q}^{t }(\Omega )\times B_{p,p}^{t-m_0 +1/p'}(\partial\Omega ),\\
\ol B_{p,q}^{t +1}(\Omega)\to \ol
B_{p,q}^{t }(\Omega )\times B_{p,q}^{t-m_0 +1/p'}(\partial\Omega ),
\endcases\tag4.7
$$
for $t>\max\{m_0,0\}-1/p'$.

$3^\circ$ The operator in {\rm (4.7)} belongs to the
Boutet de Monvel calculus; hereby the Fredholm solvability holds if
and only  if (in addition to the invertibility of the interior symbol) the
boundary symbol operator is bijective at each $(x',\xi ')\in
T^*(\partial\Omega )\setminus 0$. This can also be formulated as a
unique solvability of the model problem for {\rm (4.4)}
at each $x'\in\partial\Omega $, $\xi '\ne 0$.

$4^\circ$ In the transition between {\rm (4.4)} and {\rm (4.5)}, $\pmatrix R'_B&K'_B\endpmatrix$ is a parametrix for {\rm (4.5)} if
and only if
$$
\pmatrix R_B&K_B\endpmatrix=\pmatrix \Lambda _+^{(1-\mu
)}e^+R'_B\Lambda ^{(\mu -m)}_{-,+}&\Lambda _+^{(1-\mu
)}e^+K'_B\endpmatrix \tag 4.8
$$
is a parametrix for {\rm (4.4)}.
\endproclaim

\demo{Proof} 
The mapping (4.3) is well-defined, since $r^+B\colon X_{p,q}^{(\mu -1)(s)}(\comega)\to \ol
X_{p,q}^{s-m_0-\operatorname{Re}\mu }(\Omega )$ by Theorem 3.2 $1^\circ$,
and $\gamma _0$ acts as in (3.3).

$1^\circ$. Let us go through the transition between (4.4) and (4.5),
as already laid out in the
formulation of the theorem. 

We have from Definition 3.1  that $u\in
X_{p,q}^{(\mu -1)(s)}(\comega)$ if and only if $ w=r^+\Lambda _+^{(\mu
-1)}u\in \ol X_{p,q}^{s-\operatorname{Re}\mu +1}(\Omega )$; here
$u=\Lambda _+^{(1-\mu )}e^+w$. Moreover, since $\Lambda
^{(\varrho )}_{-,+}\colon \ol X_{p,q}^t(\Omega )\simto  \ol
X_{p,q}^{t-\operatorname{Re}\varrho }(\Omega )$ for all $\varrho $ and $t$,  $f\in \ol X_{p,q}^{s-\operatorname{Re}m
}(\Omega )$ if and only if  $g=\Lambda ^{(\mu -m)}_{-,+}f\in \ol X_{p,q}^{s-\operatorname{Re}\mu
}(\Omega )$. Hence the first equation in (4.4) carries over to
$$
\Lambda ^{(\mu -m)}_{-,+}r^+P\Lambda _+^{(1-\mu )}e^+w=g.
$$
Here $\Lambda ^{(\mu -m)}_{-,+}r^+P\Lambda _+^{(1-\mu )}e^+w$ can be
simplified to $r^+\Lambda ^{(\mu -m)}_{-}P\Lambda _+^{(1-\mu )}e^+w=P'_+w$, as
accounted for in the proof of Th.\ 4.4 in \cite{G15a} in a similar
situation.
The boundary condition in (4.4) carries over to that in (4.5) since
$B'_+w=r^+B\Lambda _+^{(1-\mu )}e^+w =r^+Bu$.

The order and type of the operators is clear from the definitions. 

$2^\circ$. Since the transition takes place by use of bijections, the
Fredholm property carries over between the two situations.

$3^\circ$. The model problem is the problem defined from the principal symbols of
the involved operators at a boundary point $x'$, in a local coordinate
system where $\Omega $ is replaced by $\rnp$ and the operator is
applied only in the $x_n$-direction for fixed $\xi '\ne 0$. The hereby
defined operator on ${\Bbb R}_+$ is in the Boutet de Monvel calculus called the boundary symbol
operator. The first statement in $3^\circ$ is just a reference to facts from
the Boutet de Monvel calculus. The second statement follows
immediately when the transition is applied on the principal symbol level. 

$4^\circ$. Finally, when 
$w=R_B'g+K_B'\psi $, then
$$
u=\Lambda _+^{(1-\mu )}e^+w=\Lambda _+^{(1-\mu )}e^+(R'_Bg+K'_B\psi )=
\Lambda _+^{(1-\mu )}e^+R'_B\Lambda ^{(\mu -m)}_{-,+}f+\Lambda _+^{(1-\mu )}e^+K'_B\psi ,
$$
showing the last statement.
\qed
\enddemo

The search for a parametrix here requires the analysis of 
model problems in Sobolev-type spaces over ${\Bbb R}_+$. It can be an advantage to reduce this question to the
boundary, where it suffices to investigate the ellipticity of a $\psi
$do (i.e., invertibility of its principal symbol), as in classical
treatments of differential and pseudodifferential problems.

\proclaim{Theorem 4.2} Consider the problem {\rm (4.3)--(4.4)} presented in Theorem {\rm
4.1}, and its transformed version {\rm (4.5)}. 

$1^\circ$ The nonhomogeneous Dirichlet
system for $P'$, $\{P'_+,\gamma _0\}$, is elliptic and has a
parametrix for $s>1/p$: $$
\pmatrix R'_D&
K'_D\endpmatrix \colon \cases \ol
F_{p,q}^{s-1 }(\Omega )\times B_{p,p}^{s-1/p}(\partial\Omega  )\to \ol F_{p,q}^{s}(\Omega ), 
\\
\ol
B_{p,q}^{s-1 }(\Omega )\times B_{p,q}^{s-1/p}(\partial\Omega  )\to \ol B_{p,q}^{s}(\Omega ).
\endcases\tag4.9
$$

$2^\circ$ Define
$$
S'_B=\gamma _0B'_+K'_D;\tag4.10
$$
a $\psi $do on $\partial\Omega $ of
order $m_0$. Then {\rm (4.3)} defines a Fredholm operator if and only
if $S'_B$ is elliptic. When it is so, and $\widetilde S'_B$ denotes a
parametrix, then $\{r^+P, \gamma _0r^+B\}$ has the parametrix $\pmatrix
R_B&K_B\endpmatrix$, where
$$
R_B= \Lambda _+^{(1-\mu )}(I-K'_D\widetilde S'_B\gamma _0B'_+)R'_D\Lambda _{-,+}^{(\mu -m)},\quad
K_B=\Lambda _+^{(1-\mu )}K'_D\widetilde S'_B.\tag 4.11
$$

\endproclaim

\demo{Proof} We first discuss the solvability of the type 0 problem 
(4.5) with $B'=I$. Set $Q_1=\Lambda _-^{(\mu -m)}P\Lambda _+^{(1-\mu )}\Lambda
_+^{(-1)}$; it is very similar to the operator $Q=\Lambda _-^{(\mu
-m)}P\Lambda _+^{(-\mu )}$ used in \cite{G15a}, Theorems 4.2 and 4.4,
being of order 0, type 0 and
having factorization index 0.  Then we can write
$$
P'=Q_1\Lambda _+^{(1)},\quad P'_+=r^+Q_1\Lambda _+^{(1)}e^+=r^+Q_1e^+r^+\Lambda _+^{(1)}e^+=Q_{1,+}\Lambda _{+,+}^{(1)},
\tag4.12
$$
where we used that $r^-\Lambda _+^{(1)}e^+$ is 0 on $\ol X_{p,q}^s(\Omega
)$ for $s>1/p$.

The operator $\Lambda _+^{(1)}$ defines an
elliptic (bijective) system for $s>1/p$, $$
\{\Lambda _{+,+}^{(1)},\gamma _0\}\colon \cases 
\ol F_{p,q}^{s}(\Omega )\simto \ol
F_{p,q}^{s-1 }(\Omega  )\times B_{p,p}^{s-1/p}(\partial\Omega  ),
\\
\ol B_{p,q}^{s}(\Omega )\simto \ol
B_{p,q}^{s-1 }(\Omega  )\times B_{p,q}^{s-1/p}(\partial\Omega  ).
\endcases\tag4.13
$$
This is shown in \cite{G90} Th.\ 5.1 for $q=2$ in the $F$-case, and extends to the
Besov-Triebel-Lizorkin spaces by the results of \cite{J96}. Composition with the operator $Q_{1,+}$
preserves this ellipticity, so $\{P'_+,\gamma _0\}$ forms an elliptic system
with regards to the mapping property
 $$
\{P'_+,\gamma _0\}
\colon \cases \ol F_{p,q}^{s}(\Omega )\to \ol
F_{p,q}^{s-1 }(\Omega  )\times B_{p,p}^{s-1/p}(\partial\Omega  ),\\
\ol B_{p,q}^{s}(\Omega )\to \ol
B_{p,q}^{s-1 }(\Omega  )\times B_{p,q}^{s-1/p}(\partial\Omega  ),
\endcases\tag4.14
$$
for $s>1/p$.
Hence there is a parametrix
$$
\pmatrix R_D'& K_D'\endpmatrix
$$
of this Dirichlet problem, 
continuous in the opposite direction of (4.14). This shows $1^\circ$.

Next, we can  discuss the general problem (4.5) by the help of this
special problem; such a discussion is standard within the Boutet de
Monvel calculus. Define $S'_B$ by (4.10),
it is a $\psi $do on $\partial\Omega $ of
order $m_0$ by the rules of calculus. If it is elliptic, it has a
parametrix that we denote $\widetilde S'_B$. 

On the principal symbol
level, the discussion takes place for exact operators; here we
denote principal symbols of the involved operators $P', B', K'_D$,
etc., by  $p', b', k'_D$, etc.
To solve the model problem (at a point $(x',\xi ')$ with $\xi '\ne
0$), with $g\in L_2(\rp)$, $\psi \in{\Bbb C}$,
$$
p'_+(x',\xi ',D_n)w(x_n)=g(x_n)\text{ on }\rp ,\quad \gamma _0b'_+(x',\xi ',D_n)w(x_n)=\psi\text{ at }x_n=0 ,\tag4.15
$$
 let $z=w-r'_Dg$,
then $z$ should satisfy
$$
p'_+z=0,\quad \gamma _0b'_+z=\psi -\gamma
_0b'_+r'_Dg \equiv \zeta .\tag4.16
$$
Assuming that $z$ satisfies the first equation, set 
$$
\gamma _0z=\varphi ;\text{ then }z=k'_D\varphi ,
$$
as the solution of the semi-homogeneous Dirichlet problem for $p'_+$.
To adapt $z$ to the second part of (4.16), we require that $\gamma _0b'_+z=\zeta $; here
$$
\gamma _0b'_+z=\gamma _0b'_+k'_D\varphi=s'_B\varphi ,
$$
when we define $s'_B$ by (4.10) on the principal symbol level; it is just
a complex number depending on $(x',\xi ')$.   The equation
$$
s'_B\varphi =\zeta \tag4.17
$$
is uniquely solvable precisely when $s'_B\ne 0$. In that case, (4.17)
is solved uniquely by $\varphi =(s'_B)^{-1}\zeta $. 

With this choice of $\varphi $, $z=k'_D\varphi $ is the unique solution of (4.16), and $w=r'_Dg+z$
is the unique solution of (4.15). The formula is in details
$$
w=r'_Dg+k'_D(s'_B)^{-1}\zeta =(I-k'_D(s'_B)^{-1}\gamma _0b'_+)r'_Dg+
k'_D(s'_B)^{-1}\psi. \tag4.18
$$

Expressed for the full operators, this shows that the problem (4.5) is
elliptic precisely when the $\psi $do $S'_B$ is so.

For the full operators, a similar construction can be carried out in
a parametrix sense, but it is perhaps simpler to test directly by compositions
that the following operator defined in analogy with (4.18):
$$
\pmatrix R'_B&K'_B\endpmatrix
=\pmatrix (I-K'_D\widetilde S'_B\gamma _0B'_+)R'_D&
K'_D\widetilde S'_B\endpmatrix\tag4.19
$$
is a parametrix for $\{P'_+, \gamma _0B'_+\}$: Since
$R'_DP'_++K'_D\gamma _0=I+\Cal R$ and $\widetilde S'_B\gamma
_0B'_+K_D=\widetilde S'_BS'_B=I+\Cal S$, with operators $\Cal R$ and
$\Cal S$ of order $-\infty $, 
$$
\aligned
\pmatrix R'_B&K'_B\endpmatrix& \pmatrix P'_+\\ \gamma _0B'_+\endpmatrix
=
(I-K'_D\widetilde S'_B\gamma _0B'_+)R'_DP'_+ +
K'_D\widetilde S'_B \gamma _0B'_+\\
&=(I-K'_D\widetilde S'_B\gamma _0B'_+)(1+\Cal R-K'_D\gamma _0) +
K'_D\widetilde S'_B \gamma _0B'_+\\
&=I-K'_D\widetilde S'_B\gamma _0B'_+ -K'_D\gamma _0+K'_D\widetilde S'_B\gamma _0B'_+K'_D\gamma _0 +
K'_D\widetilde S'_B \gamma _0B'_++\Cal R_1\\
&= I+\Cal R_2,
\endaligned\tag4.20
$$
with operators $\Cal R_1$ and $\Cal R_2$ of order $-\infty $. The
composition in the opposite order is similarly checked.

All this takes place in the Boutet de Monvel calculus. For our
original problem we now find the parametrix as in (4.11), by the
transition described in Theorem 4.1.
\qed
\enddemo

The order assumption on $B$ was made for the sake of arriving at
operators to which the Boutet de Monvel calculus applies. We think
that $m_0$ could be allowed to be noninteger, with some more
effort, drawing on results from Grubb and H\"ormander \cite{GH90}.

The treatment can be extended to problems with vector-valued boundary conditions
$\gamma _0r^+B$, when we also involve $\varrho _{\mu ,M}$ for $M> 1$,
cf.\ (3.8).

\subhead 4.2 The Neumann boundary operator $\gamma _{\mu _0-1,1}$
\endsubhead

For ease of comparison to \cite{G15a} we denote the $\mu $ used above by
$\mu _0$ here.

The boundary conditions with $B$ of noninteger order $m_0+\mu _0$ are
generally nonlocal, since $B$ is so. 
But there do exist local boundary conditions too. 
For example, the Dirichlet-type operator $\gamma _{\mu_0 -1,0}$ is local,
cf.\ (2.23). So are the systems  (cf.\ (3.8)) $\varrho _{\mu _0-M, M}=\{\gamma _{\mu _0
-M,0},\dots,\gamma _{\mu _0-M,M-1}\}$, which
also define  Fredholm operators together with $r^+P$, cf.\
Theorem 3.2 $3^\circ$.
Note that $\{r^+P,\varrho _{\mu _0-M,M}\}$
operates from a larger space $X_{p,q}^{(\mu _0-M)(s)}(\comega)$ than
$X_{p,q}^{(\mu _0-1)(s)}(\comega)$ when $M>1$.

What we shall show now is that one can impose a  {\it higher-order local}
boundary condition defined on $X_{p,q}^{(\mu _0-1)(s)}(\comega)$ {\it itself}, leading to a
meaningful boundary value problem with Fredholm solvability under a
reasonable ellipticity condition.

Here we treat the Neumann-type condition $\gamma _{\mu _0-1,1}u=\psi $,
recalling from \cite{G15a} (5.3)ff.\  that
$$
\gamma _{\mu _0-1,1}u=\Gamma (\mu _0+1)\gamma _0(\partial_n(d(x)^{1-\mu _0}u)).\tag4.21
$$
By application of (3.8) with $M=2$, $\mu =\mu _0-1$,
$$
\gamma _{\mu _0-1,1}=\gamma _{\mu ,M-1}\colon 
\cases
F_{p,q}^{( \mu _0
-1)(s)}(\comega)\to B_{p,p}^{s-\operatorname{Re}\mu
_0-1/p}(\partial\Omega ),\\
B_{p,q}^{( \mu _0
-1)(s)}(\comega)\to B_{p,q}^{s-\operatorname{Re}\mu
_0-1/p}(\partial\Omega ),
\endcases\tag4.22
$$
 is well-defined for $s>\operatorname{Re}\mu
 +M-1/p'=\operatorname{Re}\mu _0+1/p$. 

The discussion of ellipticity takes place in local coordinates, so let
us now assume that we are in a localized situation where $P$ is given on ${\Bbb R}^n$, globally estimated,
elliptic of order $m$ and of type $\mu _0$ and with factorization
index $\mu _0$ relative to the subset $\rnp$, as in \cite{G15a}, Th.\ 6.5.

For $\rnp$ we can express $\gamma _{\mu _0-1,1}$ in terms of auxiliary
operators by
$$
\gamma _{\mu _0-1,1}u=\gamma _0\partial_n\Xi _+^{\mu _0 -1}u-(\mu _0-1)[D']\gamma _0\Xi _+^{\mu _0-1}u,\tag4.23
$$
see the calculations after Cor.\ 5.3 in \cite{G15a}.
(In the manifold situation there is a certain freedom in choosing
$d(x)$ and $\partial_n$, so we are tacitly assuming
that a choice has been made that carries over to $d(x)=x_n$,
$\partial_n=\partial/\partial x_n$ in the localization.)

There is an obstacle to applying the results of Section 4.1 to this,
namely that $\Xi _+^{\mu _0-1}$ is not truly a $\psi $do! This is a
difficult fact
that has been observed throughout the development of the theory. However, in
connection with boundary conditions, operators like $\Xi _+^\mu $
work to some extent like the truly pseudodifferential operators $\Lambda _+^\mu $. It is for
this reason that we gave two versions of the operator $K_D$ in (4.2)ff.,
stemming from \cite{G15a} Th.\ 6.5 in which Lemma 6.6 there was used.

\proclaim{Theorem 4.3}  Let $P$ be given on ${\Bbb R}^n$, globally estimated,
elliptic of order $m$ and of type $\mu _0$ and with factorization
index $\mu _0$ relative to the subset $\rnp$, and let $\pmatrix
R_D&K_D\endpmatrix$ be a parametrix of the nonhomogeneous Dirichlet problem, as
recalled in {\rm (4.2)ff.}, with $K_D=\Xi _+^{1-\mu _0}e^+K'$ for a certain
Poisson operator $K'$ of order $0$.

Consider the Neumann-type problem 
$$
r^+Pu=f,\quad  \gamma _{\mu _0-1,1}u=\psi, \tag4.24
$$
where 
$$
\{r^+P,\gamma _{\mu _0-1,1}\}\colon 
\cases
F_{p,q}^{(\mu _0-1)(s)}(\crnp)\to \ol
F_{p,q}^{s-\operatorname{Re}m}(\rnp)\times
B_{p,p}^{s-\operatorname{Re}\mu _0-1/p}({\Bbb R}^{n-1}),\\
B_{p,q}^{(\mu _0-1)(s)}(\crnp)\to \ol
B_{p,q}^{s-\operatorname{Re}m}(\rnp)\times
B_{p,q}^{s-\operatorname{Re}\mu _0-1/p}({\Bbb R}^{n-1}),
\endcases\tag 4.25
$$
for $s>\mu _0+1/p$.

$1^\circ$ The operator 
$$
S_N=\gamma _{\mu _0-1,1}K_D\tag4.26
$$
equals $(\gamma _0\partial_n-(\mu _0-1)[D']\gamma _0)K'$ and is a $\psi $do on ${\Bbb R}^{n-1}$ of order $1$. 

$2^\circ$ If $S_N$ is
elliptic, then, with a parametrix of $S_N$ denoted $\widetilde S_N$,
there is the  parametrix for   $\{r^+P,\gamma _{\mu _0-1,1}\}$: 
$$
\pmatrix R_N& K_N\endpmatrix =\pmatrix (I-K_D\widetilde S_N\gamma _{\mu _0-1,1})R_D
&K_D\widetilde S_N\endpmatrix .\tag 4.27
$$
 
$3^\circ$ Ellipticity holds in particular when the principal symbol of
$P$ equals $c(x)|\xi |^{2\mu _0}$, with $\operatorname{Re}\mu
_0>0$, $c(x)\ne 0$.

\endproclaim
 
\demo{Proof} $1^\circ$. By the formulas for $\gamma _{\mu _0-1,1}$ and $K_D$, 
$$
S_N=\gamma _{\mu _0-1,1}K_D=(\gamma _0\partial_n-(\mu _0-1)[D']\gamma
_0)\Xi _+^{\mu _0-1}\Xi _+^{1-\mu _0}K' =(\gamma _0\partial_n-(\mu _0-1)[D']\gamma _0)K',
$$
and it follows from the rules of the
Boutet de Monvel calculus that this is a $\psi $do on ${\Bbb R}^{n-1}$
of order 1.

$2^\circ$. In the elliptic case, one checks that (4.27) is  a parametrix 
by  calculations as in Theorem
4.2.

$3^\circ$. In this case, the model problem for $\{r^+P,\gamma _{\mu
_0-1,1}\}$ can be reduced to that for $\{r^+(1-\Delta )^{\mu _0},\gamma
_{\mu _0-1,1}\}$. For the latter, we have shown unique solvability
in Theorem A.2 and Remark A.3 in the appendix.\qed
\enddemo

\example {Remark 4.4}
The operator $S_N$ is in fact the {\it Dirichlet-to-Neumann} operator for
$P$, sending the Dirichlet data over into the Neumann data for
solutions of $r^+Pu=0$ in an approximate sense (modulo operators of
order $-\infty $). From the calculations in the appendix we see that
its principal symbol equals $-\mu _0|\xi '|$, when $P$ is principally
equal to $(-\Delta )^{\mu _0}$, $\operatorname{Re}\mu _0>0$. 
\endexample

\subhead 4.3 Systems, further perspectives \endsubhead

The factorization property used above will not in general hold for  systems ($N\times
N$-matrices) in a convenient way  with smooth dependence on $\xi '$,
even if every element of the matrix has a factorization. But with the
$\mu $-transmission property we can establish an
extremely useful connection to systems in the Boutet
de Monvel calculus:

\proclaim{Proposition 4.5} 
Let $N$ be an integer $\ge 1$, and let $P$ be an $N\times N$-system,
 $P=(P_{jk})_{j,k=1,\dots, N}$, of  classical $\psi $do's $P_{jk} $ of order $m\in{\Bbb C}$
on $\Omega _1$ and
of type $\mu \in{\Bbb C}$ relative to 
 $\Omega
 $. Let $\mu _0\in \mu +{\Bbb Z}$. Then the operator 
$$
Q=\Lambda _-^{(\mu _0-m)}P\Lambda _+^{(-\mu _0)},\tag4.28
$$
is
of order and type $0$, and hence belongs to the Boutet de Monvel
calculus. \endproclaim

\demo{Proof} The factors $\Lambda _-^{(\mu _0-m)}$ and $\Lambda
_+^{(-\mu _0)}$ should be understood as diagonal matrices with $\Lambda _-^{(\mu _0-m)}$ resp.\ $\Lambda
_+^{(-\mu _0)}$ in the diagonal. When they are composed with $P$, they
act on each entry defining an operator of order and type 0 by the
symbol composition rules. \qed 
\enddemo

This will allow for a general application of the Boutet de Monvel
theory in the discussion of boundary value problems. Leaving the most
general case for
future works, we shall in the present paper just draw conclusions for
systems where the operator (4.28) defines a system $Q_+$ that is in
itself elliptic.
Let us give a name to such cases, where the present considerations will apply
without further efforts:

\example {Definition 4.6} Let $N$ be an integer $\ge 1$, and let $P$ be an elliptic $N\times N$-system,
 $P=(P_{jk})_{j,k=1,\dots, N}$, of  classical $\psi $do's $P_{jk} $ of order $m\in{\Bbb C}$
on $\Omega _1$ and
of type $\mu \in{\Bbb C}$ relative to 
 $\Omega
 $. Let $\mu _0\in \mu +{\Bbb Z}$. Then $P$ is said to be {\it $\mu
 _0$-reducible}, when the operator $Q$ defined in (4.28)
of order and type $0$, has the property that $Q_+$ {\it is
elliptic in the Boutet de Monvel calculus} (without auxiliary boundary operators). 
\endexample

The condition in the definition means that in local coordinates at the
boundary, the model operator
$q_0(x',0,\xi ',D_n)_+$ is bijective in $L_2(\rp)^N$.
It holds for $N=1$ for the operators with factorization index $\mu _0$, as
accounted for in the proof of \cite{G15a} Th.\ 4.4. Another important
case is where the operator $P$ (scalar or a system) is strongly
elliptic, as observed in  \cite{E81}, Ex.\ 17.1. 

\proclaim{Lemma 4.7} Let $N\ge 1$, and $P$ be of order $m\in\rp$ on
$\Omega _1$ and of type $\mu _0=m/2$ relative to $\Omega $.  
If $P$ is strongly
elliptic, i.e., satisfies in local coordinates (with $c>0$):
 $$
\operatorname{Re}(p_0(x,\xi )v,v)\ge c
|\xi |^{m}|v|^2,\text{ for all
}\xi \in{\Bbb R}^n, v\in{\Bbb C}^N,
$$
then $P$ is $\mu _0$-reducible. 
\endproclaim

\demo{Proof} Here $Q$ equals $\Lambda _-^{(-m/2)}P\Lambda
_+^{(-m/2)}$. This is  
strongly elliptic of order 0, because the principal symbols of $\Lambda
_-^{(-m/2)}$ and $\Lambda _+^{(-m/2)}$ are conjugates and
homogeneous elliptic of order $-m/2$:
 $$
\operatorname{Re}(q_0(x,\xi )v,v)=
\operatorname{Re}(p_0(x,\xi )\lambda _{+,0}^{-m/2}(\xi )v,\lambda
_{+,0}^{-m/2}(\xi )v)\ge c|\xi |^m|\lambda _{+,0}^{-m/2}(\xi )v|^2
\ge c'|v|^2,
$$
 for all
$\xi \in{\Bbb R}^n$, $ v\in{\Bbb C}^N$,
in local coordinates.
Thus at each $x'\in\partial\Omega $, $\xi '\ne 0$, the model
operator $q_0(x',0,\xi ',D_n)$ on ${\Bbb R}$ satisfies
$$
\operatorname{Re}(q_0u,u)\ge C\|u\|^2 _{L_2({\Bbb R})^N}, \text{ for
}u\in L_2({\Bbb R})^N,
$$
as seen by Fourier transformation in $\xi _n$. In particular, the
restriction of  $r^+q_0$ to $C_0^\infty (\rp)^N$
satisfies the inequality, and the inequality
extends to its closure, $r^+q_0e^+$ defined on \linebreak$L_2(\rp)^N$, which is
therefore injective. Similar
considerations hold for the adjoint, so indeed, 
 $q_0(x',0,\xi
',D_n)_+$ is bijective in $L_2(\rp)^N$. \qed

\enddemo

\proclaim{Theorem 4.8} Let $P$ be an elliptic $N\times N$ system,
$P=(P_{jk})_{j,k=1,\dots, N}$, of  classical $\psi $do's $P_{jk} $ of order $m\in{\Bbb C}$
on $\Omega _1$ and of type $\mu _0\in{\Bbb C}$ relative to $\Omega $.

Define $Q$ by {\rm(4.28)} and assume that $P$ is $\mu _0$-reducible.
Then we have:
 
$1^\circ$ Let $s>\operatorname{Re}\mu _0-1/p'$. If $u\in \dot X_{p,q}^\sigma (\comega)^N $ for
some $\sigma >\operatorname{Re}\mu _0-1/p'$ and $r^+ Pu\in \ol X_{p,q}^{s-\operatorname{Re}m}(\Omega )^N$,
then $u\in X_{p,q}^{\mu _0(s)}(\comega)^N$. 
The mapping 
$$
r^+P\colon X_{p,q}^{\mu _0(s)}(\comega )^N\to \overline
X_{p,q}^{s-\operatorname{Re}m}(\Omega )^N\tag 4.29
$$
is Fredholm, and has the parametrix 
$$
R=\Lambda ^{(-\mu
 _0)}_+e^+  \widetilde {Q_+}\Lambda ^{(\mu _0-m)}_{-,+}\colon 
\ol X_{p,q}^{s-\operatorname{Re}m}(\Omega  )^N
\to X_{p,q}^{\mu _0(s)}(\comega )^N ,  \tag4.30
$$ 
where $\widetilde { Q_+}$ is a parametrix of $Q_+$; it has the
structure $\widetilde Q_++G$ with a singular Green operator of order
and class $0$.

$2^\circ$ In particular, if $r^+ Pu\in C^\infty (\comega)^N$, then  $u\in \E_{\mu
_0}(\comega)^N$, and the mapping 
$$
r^+P\colon \E_{\mu _0}(\comega)^N\to C^\infty (\comega)^N\tag4.31
$$
is Fredholm.

$3^\circ$ Moreover, let $\mu =\mu _0-M$ for a positive integer $M$. Then when $s>\operatorname{Re}\mu
_0-1/p'$,  $\{r^+P,\varrho _{\mu ,M}\}$ 
defines a Fredholm operator
$$
\{r^+P,\varrho _{\mu ,M}\}\colon 
\cases
F_{p,q}^{\mu (s)}(\comega)^N\to \ol
F_{p,q}^{s-\operatorname{Re}m}(\Omega )^N\times \prod _{0\le
j<M}B_{p,p}^{s-\operatorname{Re}\mu -j-1/p}(\partial\Omega )^N,\\
B_{p,q}^{\mu (s)}(\comega)^N\to \ol
B_{p,q}^{s-\operatorname{Re}m}(\Omega )^N\times \prod _{0\le
j<M}B_{p,q}^{s-\operatorname{Re}\mu -j-1/p}(\partial\Omega )^N.
\endcases\tag4.32
$$
\endproclaim

\demo{Proof}
The proof goes as in \cite{G15a} Theorems 4.4 and 6.1:

$1^\circ$. We replace the equation
$$
r^+ Pu=f\in \ol X_{p,q}^{s-\operatorname{Re}m}(\Omega )^N,\tag 4.33
$$
by composition to the left with $\Lambda ^{(\mu _0-m)}_{-,+}$, by the equivalent
problem
$$
\Lambda ^{(\mu _0-m)}_{-,+}r^+ Pu=g,\text{ where }g=\Lambda ^{(\mu _0-m)}_{-,+}f\in \ol X_{p,q}^{s-\operatorname{Re}\mu _0}(\Omega )^N,\tag 4.34
$$
using the homeomorphism properties of $\Lambda ^{(\mu _0-m)}_{-,+}$,
applied to vectors. Here
$f=\Lambda ^{(m-\mu _0)}_{-,+}g$. Moreover,
cf.\ Remark 1.1 in \cite{G15a},
$$
\Lambda ^{(\mu _0-m)}_{-,+}r^+ Pu=r^+\Lambda ^{(\mu _0-m)}_{-} Pu.
$$
Next, we set
 $v=r^+\Lambda ^{(\mu _0)}_+u$; then $u=\Lambda ^{(-\mu
 _0)}_+e^+v$, and equation (4.33) becomes 
$$
Q_+v=g; \quad g\text{ given in }\ol X_{p,q}^{s-\operatorname{Re}\mu _0}(\Omega ),\tag4.35
$$
where $Q$ is defined by (4.28).

The properties of $P$ imply that $Q$ is elliptic of order 0 and type
0, hence belongs to the Boutet de
Monvel calculus. The rest of the argumentation takes place within that
calculus.
By our assumption, 
$Q_+=r^+Qe^+$ defines an elliptic boundary problem (without auxiliary trace or
Poisson operators) there, and $Q_+$ is continuous in  $\ol
X_{p,q}^t(\Omega )$ for $t>-1/p'$.  By the ellipticity,  $Q_+$ has a
parametrix $\widetilde {Q_+}$, continuous in the opposite direction,
and with the mentioned structure. Since $v\in \dot X_{p,q}^{-1/p'+0}(\comega)$ by
hypothesis,  solutions of $Q_+v=g$ with $g\in \ol X_{p,q}^t(\Omega )$
for some $t>-1/p'$ are in $\ol X_{p,q}^t(\Omega )$. Moreover,
$$
Q_+\colon \ol X_{p,q}^t(\Omega )\to \ol X_{p,q}^t(\Omega )\text{ is Fredholm for all }t>-1/p'.
$$
When carried back to the original functions, this shows $1^\circ$.

$2^\circ$ follows by letting $s\to\infty $, using that $\bigcap _{s}X_{p,q}^{\mu (s)}(\comega)^N=\E_\mu
(\comega)^N$.

For $3^\circ$, we use that the mapping $\varrho _{\mu ,M}$ in (3.8) extends immediately to
vector-valued functions:
$$
\varrho _{\mu ,M}\colon \cases
F_{p,q}^{\mu (s)}(\comega)^N
\to  \prod _{0\le j<M}B_{p,p}^{s-\operatorname{Re}\mu
-j-1/p}(\partial\Omega )^N,
\\
B_{p,q}^{\mu (s)}(\comega)^N
\to  \prod _{0\le j<M}B_{p,q}^{s-\operatorname{Re}\mu
-j-1/p}(\partial\Omega )^N,
\endcases\tag4.36
$$
when $s>\operatorname{Re}\mu _0 -1/p'$; surjective with nullspace $X_{p,q}^{\mu _0(s)}(\comega)^N$ (recall $\mu =\mu _0-M$). When we adjoin this
mapping to (4.29), we obtain (4.32).
\qed  
\enddemo

One of the things we obtain here is that results from
 \cite{E81} (extended to $L_p$ in \cite{S95, CD01}), on solvability for $s$ in an interval of length 1 around
 $\operatorname{Re}\mu _0$, are lifted to regularity and Fredholm
 properties for all larger $s$, with exact
information on the domain,  also in 
general scales of function spaces. Moreover, our theorem is obtained via
a systematic variable-coefficient calculus, whereas the results in
\cite{E81} are derived from constant-coefficient considerations by ad
hoc perturbation methods in $L_2$-Sobolev spaces.

Also the results on other boundary conditions in the present paper 
extend to suitable systems. One can moreover extend the results to operators in
vector bundles (since they are locally matrix-formed).

\medskip
The Boutet de Monvel theory is not an
easy theory (as the elaborate presentations \cite{B71, RS82,
G84, G90, G96, S01, G09} in the literature show), but
one could have feared that a theory for the more general $\mu
$-transmission operators and their boundary problems would be a step up
in difficulties. Fortunately, as we have seen, many of the issues can be dealt
with by reductions using the special operators $\Lambda _\pm^{(\mu )}$,
to cases where the type 0 theory applies.

There is currently also an interest for problems with less smooth
symbols. In this connection we mention
that there do exist pseudodifferential theories for such problems,
also with boundary conditions, cf.\ Abels
\cite{A05} and \cite{G14} and their references. One finds that a lack of
smoothness in the $x$-variable narrows down the interval of parameters
$s$  where one has good solvability properties, and
compositions are delicate. --- It is also possible to work under
limitations on the number of standard estimates in $\xi $.

\head Appendix. Calculations in an explicit example \endhead

Pseudodifferential methods are a refinement of the application of the Fourier
transform, making it useful even for variable-coefficient partial
differential operators, and allowing generalizations to e.g.\ operators
of noninteger order. But to explain some basic mechanisms it may be
useful to consider a simple ``constant-coefficient'' case, where
explicit elementary calculations can be made, not requiring intricate
composition rules.
 This is the case for
$(1-\Delta )^{a}$ ($a>0$) on $\rnp$, where everything
can be worked out by hand in exact detail (in the spirit of the elementary Ch.\
9 of [G09]). We here restrict the attention to $H^s_p$-spaces.

The symbol of $(1-\Delta )^a$ is factorized as
$$
(\ang{\xi '}^2+\xi _n^2)^a=(\ang{\xi '}-i\xi _n)^a(\ang{\xi '}+i\xi _n)^a.\tag A.1
$$

Now we shall use the definitions of simple order-reducing
operators $\Xi ^t_\pm$ and Poisson operators $K_j$ from [G15a] with $\ang{\xi '}$ instead
of $[\xi ']$, because they fit particularly well with the factors in
(A.1).
We shall often abbreviate $\ang{\xi '}$ to $\sigma $.

The {\it homogeneous Dirichlet problem}
$$
r^+(1-\Delta  )^au=f, \quad f \text{ given in }\ol H^{s-2a}_p(\rnp),\tag A.2
$$
$s>a-1/p'$, has a unique solution $u$ in $ \dot
H_p^{a-1/p'+0}(\crnp)$ determined as follows:

With $\Xi ^t_{\pm}=\operatorname{OP}\bigl((\ang{\xi '}+i\xi _n)^t\bigr)$,
we have that
 $(1-\Delta )^a=\Xi ^a_-\Xi ^a_+$ on ${\Bbb R}^n$. Let $v=r^+\Xi ^a_+u$; it
is in $\ol H_p^{\,-1/p'+0}(\rnp)= \dot H_p^{-1/p'+0}(\crnp)
$, and $u=\Xi _+^{-a}e^+v$. Then (A.2) is turned into
$$
r^+\Xi ^a_-e^+v=f.\tag A.3
$$
Here $r^+\Xi ^a_-e^+=\Xi ^a_{-,+}$ is known to map $\ol H^{t}_p(\rnp)$
homeomorphically onto $\ol H^{t-a}_p(\rnp)$ for all $t\in{\Bbb R}$,
with inverse $\Xi ^{-a}_{-,+}$. (Cf.\ e.g.\ \cite{G15a} Sect.\ 1.) In particular, with $f$ given in $\ol
H^{s-2a}_p(\rnp)$, (A.3) has the unique solution $v=\Xi
^{-a}_{-,+}f\in \ol H^{s-a}_p(\rnp)$. Then (A.2) has the unique solution
$$
u=\Xi _+^{-a}e^+\Xi ^{-a}_{-,+}f\equiv R_Df,\tag A.4
$$
and it belongs to $H_p^{a(s)}(\crnp)$ by the definition of that
space. Thus the solution operator for (A.2) is $R_D=\Xi _+^{-a}e^+\Xi
^{-a}_{-,+}$. (This is a simple variant of the proof of [G15a] Th.\ 4.4.)

Next, we go to the larger space $H_p^{(a-1)(s)}(\crnp)$, still
assuming $s>a-1/p'$, where we study the {\it nonhomogeneous Dirichlet
problem}. By [G15a] Th.\ 5.1
with $\mu =a-1$ and $M=1$, we have a mapping $\gamma _{a-1,0}$, acting as
$$
\gamma _{a-1,0}\colon u\mapsto \Gamma (a)\gamma _0(x_n^{1-a} u),
$$
also equal to $\gamma _0\Xi ^{a-1}_+u$, and sending $H_p^{(a-1)(s)}(\crnp)$ onto
$B_p^{s-a+1-1/p}({\Bbb R}^{n-1})$ with kernel
$H_p^{a(s)}(\crnp)$. Together with $(1-\Delta )^a$ it therefore
defines a homeomorphism for $s>a-1/p'$:
$$
\{r^+(1-\Delta )^a,\gamma _{a-1,0}\}\colon H_p^{(a-1)(s)}(\crnp)\to \ol H_p^{s-2a}(\rnp)\times B_p^{s-a+1-1/p}({\Bbb R}^{n-1}).\tag A.5
$$
It represents  the problem
$$
r^+(1-\Delta )^au=f, \quad  \gamma _{a-1,0}u=\varphi ,\tag A.6
$$
that we regard as the nonhomogeneous Dirichlet problem for
$(1-\Delta )^a$.
The solution operator in the case $\varphi =0$ is clearly $R_D$ defined above,
since the kernel of $\gamma _{a-1,0}$ is $H_p^{a(s)}(\crnp)$.

Also the solution operator for the problem (A.6) with $f=0$
can be found explicitly: 

On the boundary symbol level we consider
the problem (recall $\sigma =\ang{\xi '}$)
$$
(\sigma -\partial _n)^a(\sigma +\partial_n)^au(x_n)=0\text{ on }\rp.\tag A.7
$$
Since $\operatorname{OP}_n((\sigma -i\xi _n)^\mu )$ preserves support in
$\crm$ for all $\mu $, $u$ must equivalently satisfy
$$
(\sigma +\partial_n)^au(x_n)=0\text{ on }\rp.\tag A.8
$$
This has the distribution solution
$$
u(x_n)=\F^{-1}_{\xi _n\to x_n}(\sigma +i\xi _n)^{-a}=
\Gamma (a)^{-1}x_n^{a-1}e^+r^+e^{-\sigma x_n}\tag A.9
$$
(cf.\ e.g.\ \cite{H83} Ex.\ 7.1.17 or [G15a] (2.5)), and the derivatives $\partial_n^ku$ are likewise solutions, since 
$$
(\sigma +i\xi _n)^{a}(i\xi _n)^k(\sigma +i\xi _n)^{-a}=(i\xi
_n)^k=\F_{x_n\to\xi _n}\delta _0^{(k)},
$$ where $\delta _0^{(k)}$ is supported in
$\{0\}$. The undifferentiated function matches our problem. Set
$$
\tilde k_{a-1,0}(x_n,\xi ')=\Gamma (a)^{-1}x_n^{a-1}e^+r^+e^{-\sigma x_n}=
\F^{-1}_{\xi _n\to x_n} (\sigma +i\xi _n)^{-a},
\tag A.10
$$
then since $\gamma _{a-1,0}\tilde k_{a-1,0}=1$, the mapping ${\Bbb C}\ni \varphi \mapsto
\varphi \cdot r^+ \tilde k_{a-1,0}$ solves the problem
$$
(\sigma +\partial_n)^au(x_n)=0\text{ on }\rp,\quad \gamma _{a-1,0}u=\varphi .\tag A.11
$$
Using the Fourier transform in
$\xi '$ also, we find that (A.6) with $f=0$ has the solution 
$$
u(x)=K_{a-1,0}\varphi \equiv\F^{-1}_{\xi '\to x'}\bigl(\tilde k_{a-1,0}(x_n,\xi ')\hat\varphi (\xi ')\bigr).\tag A.12
$$
It can be denoted $\operatorname{OPK}(\tilde k_{a-1,0})\varphi $, by a
generalization 
of the notation from the Boutet de Monvel calculus. We moreover define
$k_{a-1,0}(\xi )=\F_{x_n\to\xi _n}\tilde k_{a-1,0}(x_n,\xi ')=(\sigma +i\xi _n)^{-a}$; $\tilde k_{a-1,0}$ and $k_{a-1,0}$ are
the symbol-kernel and symbol of $K_{a-1,0}$, respectively.

Note that
$$
\aligned
k_{a-1,0}(\xi ',\xi _n)&= (\ang{\xi '}+i\xi _n)^{-a}= (\ang{\xi '} +i\xi _n)^{1-a}(\ang{\xi '}+i\xi _n)^{-1},\text{
hence }\\
K_{a-1,0}&= \Xi
^{1-a}_+K_0,
\endaligned
\tag A.13 
$$
where $K_0=\operatorname{OPK}((\ang{\xi '}+i\xi _n)^{-1})$ is
{\it the Poisson operator for the Dirichlet problem for} $1-\Delta $,
$$
K_0\varphi =\F^{-1}_{\xi \to x}((\ang{\xi '}+i\xi _n)^{-1}\hat\varphi (\xi ')),
$$
(cf.\ e.g.\ \cite{G09}, Ch.\ 9). $K_0$ is usually considered as
mapping into a space over $\rnp$, and it is well-known that $K_0\colon B_p^{t-1/p}({\Bbb
R}^{n-1})\to  \ol H_p^{t}(\rnp)$ for all $t\in{\Bbb R}$. However, the
above formula shows that it in fact maps into distributions on ${\Bbb
R}^n$ supported in $\crnp$, so we can, with a slight abuse of notation, identify $K_0$ with $e^+K_0$,
mapping into $e^+\ol H_p^{t}(\rnp)$, and conclude that 
$$
K_{a-1,0}
\colon B_p^{s-a+1-1/p}({\Bbb
R}^{n-1})\to  H_p^{(a-1)(s)}(\crnp), \text{ for all }s\in{\Bbb R}.\tag A.14
$$

We have shown:

\proclaim{Theorem A.1} Let $a>0$. The nonhomogeneous Dirichlet problem {\rm (A.6)} for
$(1-\Delta )^a$ on $\rnp$ is uniquely solvable, 
in that the operator
{\rm (A.5)} for $s>a-1/p'$ has the inverse
$$
\pmatrix r^+(1-\Delta )^a\\ \gamma _{a-1,0}\endpmatrix^{-1}=\pmatrix R_D&K_{a-1,0}\endpmatrix,\tag A.15
$$
where $R_D$ and $K_{a-1,0}$ are defined in {\rm (A.4)} and {\rm (A.12)}. 
\endproclaim

Third, we consider the boundary problem
$$
r^+(1-\Delta )^au=f, \quad  \gamma _{a-1,1}u=\psi , \tag A.16
$$
that we shall view as a {\it nonhomogeneous Neumann problem} for $(1-\Delta )^a$. We here assume
$s>(a-1)+2-1/p'=a+1/p$, to use the construction in [G15a] Th.\ 5.1 with
$\mu =a-1$, $M=2$.  Recall from [G15a] (5.3)ff., that $\gamma _{a-1,1}$
acts as
$$
\gamma _{a-1,1}\colon u\mapsto \Gamma (a+1)\gamma _0(\partial_n(x_n^{1-a}u)).\tag A.17
$$
Moreover, we can infer from [G15a], the text after Cor.\ 5.3 (with $[\xi ']$ replaced by $\ang{\xi
'}$), that
$$
\gamma _{a-1,1}u=\gamma _0\partial_n\Xi ^{a-1}_+u-(a-1) \,\ang{D'}\,\gamma _{a-1,0}u,
$$
for $u\in H_p^{(a-1)(s)}(\crnp)$ with $s>a+1/p$. Then for a null
solution $z$ written in the form $z=K _{a-1,0}\varphi
=\Xi ^{1-a}_+K_0\varphi $ (recall (A.13)), we have since $\gamma _0\partial_n K_0=-\ang{D'}$, 
$$
\gamma _{a-1,1}z=\gamma _0\partial_n\Xi
^{a-1}_+z-(a-1)\,\ang{D'}\,\gamma _{a-1,0}z
=\gamma _0\partial_nK_0\varphi 
 -(a-1) \ang{D'}\varphi =-a\,\ang{D'}\varphi .
$$
Hence in order for $z$ to solve (A.16) with $f=0$, $\varphi $ must satisfy
$$
\psi =-a\,\ang{D'}\varphi .
$$
Since $a\ne 0$, the coefficient $-a\ang{D'}$ is an elliptic invertible
$\psi $do, so  (A.16) with $f=0$
is uniquely solvable with solution 
$$
z=K_N\psi ,\text{ where }K_N=-K_{a-1,0}a^{-1}\ang{D'}^{-1}=-\Xi _+^{1-a}K_0a^{-1}\ang{D'}^{-1}
 .\tag A.18
$$

To solve (A.16) with a given $f\ne 0$, and $\psi =0$, we let $v=R_Df$ and
reduce to the problem for $z=u-v$:
$$
r^+(1-\Delta )^a(u-v)=0,\quad \gamma _{a-1,1}(u-v)=-\gamma _{a-1,1}R_Df.
$$
This has the unique solution
$$
u-v=-K_N\gamma _{a-1,1}R_Df; \text{ hence }u=R_Df-K_N\gamma _{a-1,1}R_Df.
$$
Altogether, we find:

\proclaim{Theorem A.2}  The Neumann problem {\rm (A.16)} for
$(1-\Delta )^a$ on $\rnp$ is uniquely solvable, in that the
operator 
$$
\{r^+(1-\Delta )^a,\gamma _{a-1,1}\}\colon H_p^{(a-1)(s)}(\crnp)\to \ol H_p^{s-2a}(\rnp)\times B_p^{s-a-1/p}({\Bbb R}^{n-1}),\tag A.19
$$
for $s>a+1/p$ is a homeomorphism, with inverse
$$
\pmatrix R_N& K_N\endpmatrix =\pmatrix (I-K_N\gamma _{a-1,1})R_D
&K_N\endpmatrix ,\tag A.20
$$
with $R_D$ and $K_N$ described in {\rm (A.4)} and {\rm (A.18)}.

\endproclaim

Note that there is here a {\it Dirichlet-to-Neumann operator} $P_{DN}$
sending the Dirichlet-type data over into Neumann-type data for solutions of $r^+(1-\Delta )^au=0$:
$$
P_{DN}=-a\ang{D'}.\tag A.21
$$

\example{Remark A.3}
We have here assumed $a$ real in order to relate to the fractional powers
of the Laplacian, but all the above goes through in the same way if
$a$ is replaced by a complex $\mu $ with $\operatorname{Re}\mu >0$;
then in Sobolev exponents and inequalitites for $s$, $a$ should be
replaced by $\operatorname{Re}\mu $.
\endexample

One can also let higher order boundary operators $\gamma _{a-1,j}$
enter in a similar way, defining single boundary conditions.

\head{Acknowledgement}\endhead We are grateful to J.\
Johnsen and X.\ Ros-Oton for useful discussions.

\enddocument

\end

\bye